\documentclass[preprint,12pt]{elsarticle}

\usepackage{amssymb}
\usepackage{amsmath}
\newcommand{\G}{\Gamma}
\newcommand{\Gs}{\Gamma^*}
\newcommand{\Gi}{\Gamma_\infty}
\newcommand{\Gis}{\Gamma_\infty^*}
\newcommand{\Gstt}{\Gamma^*(\tilde t)}

\newcommand{\D}{{\bf D}}
\newcommand{\J}{{\bf J}}
\newcommand{\F}{{\bf F}}
\newcommand{\al}{{\alpha}}
\newcommand{\be}{{\beta}}
\newcommand{\bt}{\tilde \beta}
\newcommand{\ttt}{\tilde t}
\newcommand{\ts}{t^*}
\newcommand{\tG}{\tilde \Gamma}
\newcommand{\tb}{\tilde \beta}
\newcommand{\tP}{\tilde \Pi}
\newcommand{\er}{\mathrm{erfc}}
\newcommand{\Ba}{{ Ba}}

\newcommand{\ds}{\displaystyle}

\journal{}

\begin{document}

\begin{frontmatter}



\title{Fractional calculus approach to models of adsorption: Barrier-diffusion control}


\author[label1,label2]{Ivan Bazhlekov
\corref{cor1}
}
\ead{i.bazhlekov@math.bas.bg}

\author[label1,label2]{Emilia Bazhlekova}
\ead{e.bazhlekova@math.bas.bg}

\cortext[cor1]{Corresponding author}
\affiliation[label1]{organization={Institute of Mathematics and Informatics, Bulgarian Academy of Sciences},
            addressline={Acad.~G.~Bonchev str., bld. 8}, 
            city={Sofia},
            postcode={1113}, 
            country={Bulgaria}}
						
	\affiliation[label2]{organization={Centre of Excellence in Informatics and Information and Communication
						Technologies},
           city={Sofia},
           country={Bulgaria}}					
%Correspondence {i.bazhlekov@math.bas.bg (I.B.); e.bazhlekova@math.bas.bg (E.B.)}
%% Abstract
\begin{abstract}
The mathematical model of surfactant adsorption under mixed barrier-diffusion control is analyzed using techniques from fractional calculus. 
The kinetic models of Henry, Langmuir, Frumkin, Volmer and van der Waals are considered. 
First, treating the Ward-Tordai integral equation as a fractional order one, the partial differential model is transformed 
into a single fractional ordinary differential equation for the adsorption. 
A transformation of the obtained equation is proposed that reduces the number of 
parameters to two dimensionless groups (at Frumkin and van der Waals models a third parameter appears). 
In the simplest case of Henry adsorption isotherm the fractional differential model depends on a single dimensionless group and an exact solution exists, represented in terms of Mittag-Leffler functions. Based on this solution, second order asymptotes (at small values of the adsorption) are derived for the other models.  
The asymptotes of the adsorption result in a higher order asymptotes for the surface pressure (surface tension). 
For small surface coverage, all considered models converge to the Henry model's predictions, making it a universal first-order approximation for the surface tension.
Next, the fractional differential model is written as an integral equation %of fractional order 
that can be considered as a generalization of the well-known Ward-Tordai equation to the case of barrier-diffusion control. 
For computer simulation of the obtained integral equation a predictor-corrector numerical method is developed and numerical results are presented and discussed.
\end{abstract}

\begin{graphicalabstract}
\begin{figure}[h]%% placement specifier
\centering%% For centre alignment of image.
\includegraphics[width=16cm]{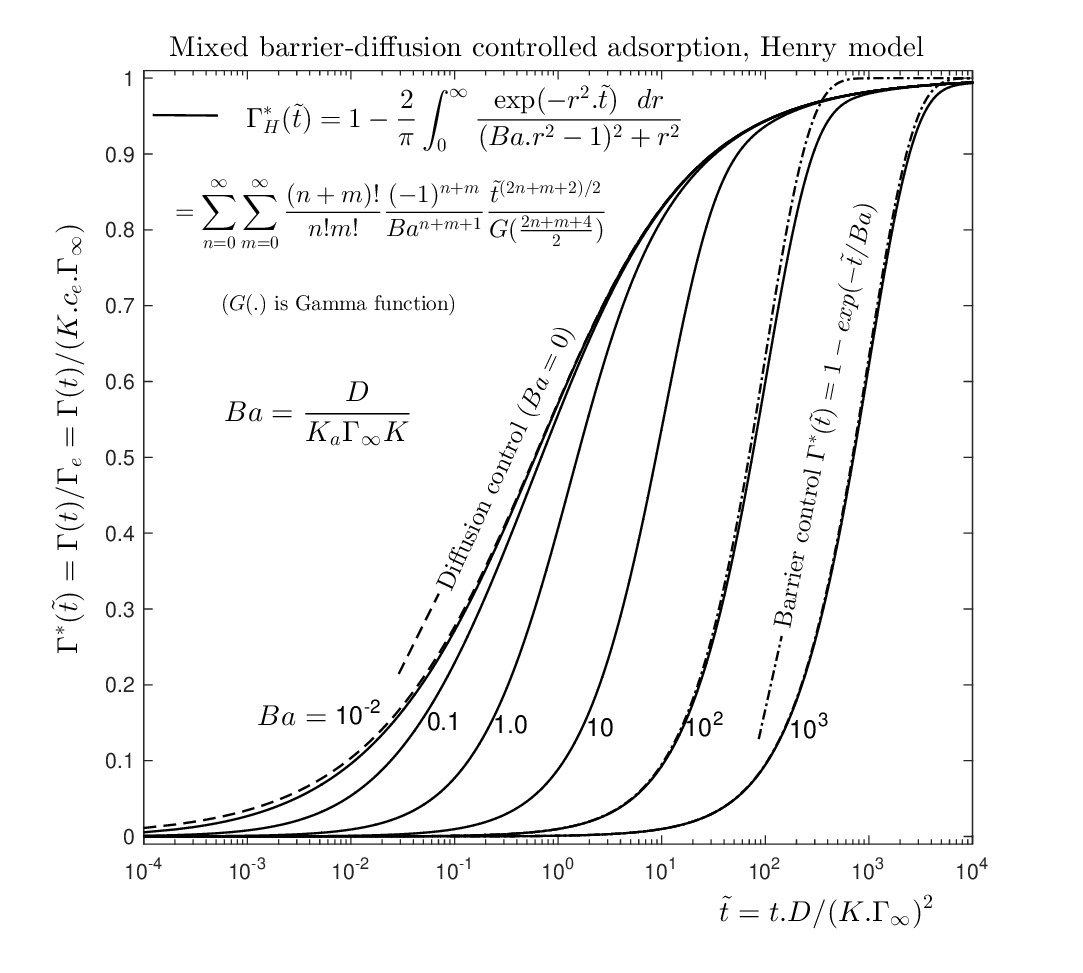}
%\includegraphics[width=12cm]{FigGrAbstr2.eps}
%%\caption{Figure Caption}\label{fig1}
\end{figure}
\begin{figure}[h]%% placement specifier
\centering%% For centre alignment of image.
\includegraphics[width=16cm]{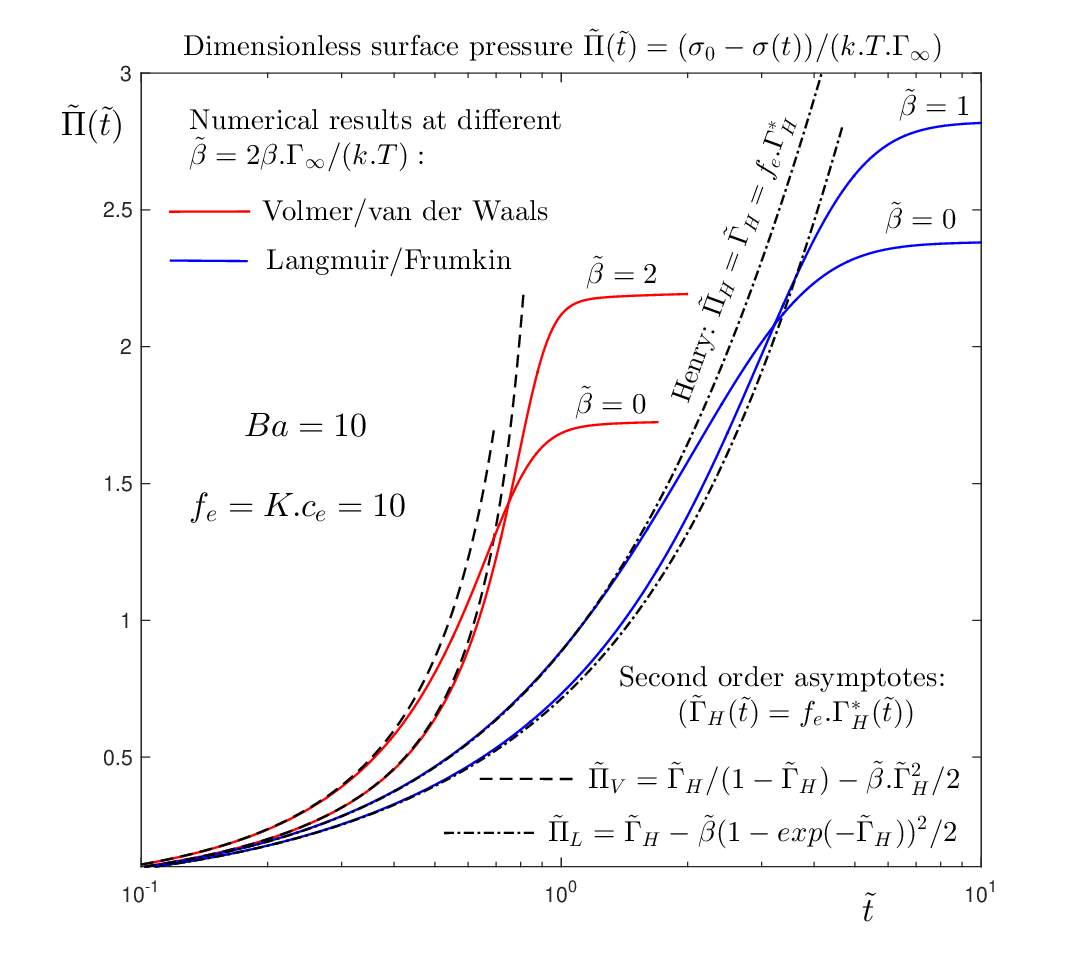}
%%\caption{Figure Caption}\label{fig1}
\end{figure}

\end{graphicalabstract}

\begin{highlights}
\item Two-term fractional ordinary differential equation for the surfactant adsorption under mixed barrier-diffusion control.
\item Generalized Ward-Tordai integral equation in the case of mixed kinetics.
\item Fractional order integral equation for the adsorption evolution.
%\item Exact analytical solutions in the case of linear Henry kinetic model.
\item Second order asymptotes for the adsorption and the surface tension at small surface coverage.
\item Common, for all considered models, first-order asymptote for the surface tension.
%\item Second order asymptotes of the adsorption at small surfactant concentrations
%are derived ( values of the adsorption, $\Gamma << \Gamma_{\infty}$)
%\item Unified, independent on the adsorption model, asymptote for the surface tension at small surface coverage
\end{highlights}

%% Keywords

\end{frontmatter}



\section {Introduction}
\label{secIn}

The presence of surfactants in a multiphase fluid system is essential in many man-made and natural processes. 
When adsorbed on an interface the surfactant not only reduces the surface tension, but also introduces important mechanical properties, 
such as surface elasticity and viscosity. 
This importance motivates the large number of experimental (see for instance \cite{CNHL2018, LTL2021, LZLWL2022, BMDSC2024}), 
theoretical \cite{A1987, ED2000, MS2015, SBPR2023}, and numerical \cite{MK1980, LFB2003, CCE2010, CFMN2015} studies of the adsorption.
For better understanding of the basic mechanisms of the process often the experimental investigations are analyzed using theoretical predictions 
\cite{DGSV2021, LPL2022, HMY2025, HIR2025}.

%In order to analyze the effect of different situation on the adsorption, different variations are studied: 
The process of surfactant adsorption is studied in different situations:
spherical interfaces are investigated in \cite{CCE2010, LHSL2022, KS2023}, 
a general expression of the adsorption in the case of expanding interface is presented in \cite{LYZM2005}, 
surfactant transport and adsorption-desorption at liquid-liquid interface are studied in \cite{LRF1997},
adsorption under the assumption of anomalous diffusion is considered in \cite{TNL2020, BB2025}. 

The pioneering work of Ward and Tordai (1946) \cite{WT1946} is the first theoretical study on the dynamics of surfactant adsorption. 
The Ward-Tordai integral equation is still a basis for most of the theoretical studies that analyze the adsorption evolution \cite{MS2015, MAF2017}.
Fifteen years later Hansen \cite{H1961}, using Laplace transform technique, derives an integro-differential equation as a model of the same problem. 
It was pointed out recently by Hristov \cite{H2016} that Ward-Tordai and Hansen equations 
are, respectively, fractional integral and differential equations of order $1/2$. 
Thus, the famous Ward-Tordai integral equation and the Hansen integro-differential equation, governing the 
diffusion-controlled adsorption of surfactant at air/liquid interface, 
are remarkable examples of fractional-order equations derived from a classical integer-order model.

Elements from the theory of fractional calculus to models of adsorption have been previously used. 
Miura and Seki \cite{MS2015} in their study of barrier-diffusion controlled adsorption use the Gr\"unwald-Letnikov formula to approximate fractional derivative of order $1/2$ and investigate the effects of anomalous subdiffusion, modeling the latter by employing Riemann-Liouville time-fractional derivative in the diffusion equation.

In the present study models of mixed barrier-diffusion controlled adsorption are analyzed using techniques from fractional calculus.
%The basis is that Ward-Tordai and Hansen equations are considered of fractional order 1/2, integral and differential, respectively.
The model of adsorption is written in Section~\ref{SecFDM} as ordinary two-term differential equation, containing a fractional derivative. 
In Section~\ref{SecFIM} the model is transformed into an integral equation %of fractional order 
for the adsorption. 
Based on existing analytical solutions to linear fractional order equations, solutions in the case of Henry adsorption isotherm are derived.
They are used in Section~\ref{SecAs} for obtaining of second order asymptotes (at small surface coverage) for the adsorption and the surface tension.
The asymptotes are validated by comparisons with results, obtained by numerical simulations. 
The numerical approach is described in Section~\ref{SecNM}.

\section{Preliminaries}
\label{secPr}

This section contains the initial mathematical model -- 
the well-known system of partial differential equations governing adsorption of surfactant under 
mixed barrier-diffusion control and the related Ward-Tordai equation. 
Basic definitions and properties of fractional calculus operators are presented as well.

\subsection{Mathematical model}
\label{subsecMM}

The mathematical model is given below in the simplest case: mass transfer of surfactant due to Fickian diffusion in a semi-infinite domain ($x>0$) 
in a quiescent bulk phase onto freshly created flat interface \cite{MAF2017, LPL2022, N2022}. The governing equation for the diffusion is
\begin{equation}
\frac{\partial c(x,t)}{\partial t}=D\frac{\partial^2 c(x,t)}{\partial x^2},\ \ \ t>0,\ x>0,
\label{DE}
\end{equation}
where $c(x,t)$ is the surfactant concentration in the bulk phase and $D$ is the diffusion coefficient. 

Initially, the surfactant is uniformly distributed with concentration $c_e$.
The initial and boundary conditions on the interface $(x=0)$, and at infinity are the following:
\begin{eqnarray}
c(x>0,t=0)&=&c_{e},\nonumber \\
c(x\to \infty,t>0)&=&c_{e},\label{IBC} \\
c(x=0,t>0)&=&c_s(t),\nonumber
\end{eqnarray}
where $c_s(t)$ is the subsurface surfactant concentration. 

The change of surfactant concentration $\G(t)$ on the interface is due to the diffusive flux from the bulk:
\begin{equation}
 \frac{d\G(t)}{dt}=D\left.\frac{\partial c(x,t)}{\partial x}\right|_{x=0},\ \ \ t>0, \ \ \ \G(0)=0.
\label{Ads}
\end{equation}

To close the mathematical model (\ref{DE}-\ref{Ads}) additional relation between $\G(t)$ and $c_s(t)$ is necessary. 
In the case of mixed barrier-diffusion controlled adsorption it represents the total rate of adsorption 
as the difference between the adsorption and desorption fluxes \cite{DVK2002, KDD2016}. 
Thus, $\G(t)$ and $c_s(t)$ are related via the differential equation:

\begin{equation}
 \frac{d\G(t)}{dt}=K_a \Phi(\G(t))\left[c_s(t) - f(\G(t))/K\right].
\label{BDCA}
\end{equation}
Here, $K_a$ is the rate constant of adsorption, the constant $K$ is the adsorption parameter 
that characterizes the surface activity of the surfactant. 
The functions $\Phi(\G)$ and $f(\G)$ describe the type of the kinetic model of adsorption. They are given in 
Table~\ref{Table1} for some of the most popular adsorption isotherms, where 
the saturation adsorption $\Gi$ represents the maximum possible surface concentration, the parameter $\be$ accounts 
for the intermolecular interaction between adsorbed molecules, $k$ is the Boltzmann constant and $T$ the absolute temperature.
The other notations used in Table~\ref{Table1} slightly differ from that in \cite{ED2000, DVK2002}: 
the adsorption constant $K_a$ is moved to eq. (\ref{BDCA}), which is done in order to simplify further notations.
The dependence of the surface pressure $\Pi(t)$ (the deviation of the dynamic surface tension $\sigma (t)$ from that of clean interface $\sigma_0$)
on the surfactant concentration $\G (t)$ is also given in Table~\ref{Table1}, which will be used further in the study. 

\begin{table}[t]
\centering%% For centre alignment of tabular.
%\begin{tabular}{|c|c|c|}
\begin{tabular}{l l l l}%% Table column specifiers
\hline
Adsorption  & $\Phi(\G(t))$ & $f(\G(t))$ & $\Pi(t)=\sigma_0-\sigma(t)$\\
isotherm & & & $= k.T.J(\G(t))$,\\
         & & & where $J(\G (t))=$\\
\hline
Henry &  $1$ & $\frac{\G(t)}{\Gi}$ & $\G(t)$\\
Langmuir  &  $1-\frac{\G(t)}{\Gi}$ & $\frac{\G(t)}{\Gi-\G(t)}$ & $-\Gi\ln(1-\frac{\G(t)}{\Gi})$\\
Frumkin  & $1-\frac{\G(t)}{\Gi}$ & $\frac{\G(t)}{\Gi-\G(t)}\exp\left(-\frac{2 \be.\G(t)}{k.T} \right)$ & $-\Gi\ln(1-\frac{\G(t)}{\Gi})-\frac{\be.\G^2(t)}{k.T} $\\
Volmer  & $1$ & $\frac{\G(t)}{\Gi-\G(t)}\exp\left(\frac{\G(t)}{\Gi-\G(t)}\right)$ & $\frac{\Gi .\G(t)}{\Gi-\G(t)}$ \\
van der Waals  & $1$ & $\frac{\G(t)}{\Gi-\G(t)}\exp\left(\frac{\G(t)}{\Gi-\G(t)}-\frac{2 \be.\G(t)}{k.T} \right)$ & $\frac{\Gi .\G(t)}{\Gi-\G(t)}-\frac{\be.\G^2(t)}{k.T}$ \\
\hline
\end{tabular}
\caption{Expressions for $\Phi(\G)$, $f(\G)$ and $J(t)$ for the considered kinetic models}\label{Table1}
\end{table}

The limiting case of diffusion-controlled adsorption, in which  
(\ref{BDCA}) reduces to a direct relation between $c_s(t)$ and $\G(t)$:
\begin{equation}
 c_s(t) = f(\G(t))/K
\label{DCA}
\end{equation}
 is widely investigated. In this case the mathematical model (\ref{DE}-\ref{Ads}) implies the following relation, derived for the first time by Ward and Tordai \cite{WT1946}:
    \begin{equation}    
\G(t) = \sqrt{\frac{D}{\pi}}\int_0^t \frac{c_e-c_s(\tau)}{(t-\tau)^{1/2}}\, d\tau,
\label{WTI}
      \end{equation}
known as Ward-Tordai equation. 
It is extensively applied for solution of the problem in the case of diffusion-controlled adsorption. 
Indeed, plugging in (\ref{WTI}) of $c_s(t)$ given by (\ref{DCA}) yields a single Volterra integral equation for $\G(t)$:
    \begin{equation}    
\G(t) = \sqrt{\frac{D}{\pi}}\int_0^t \frac{c_e-f(\G(\tau))/K}{(t-\tau)^{1/2}}\, d\tau.
\label{WTFI}
      \end{equation}

Another form of the solution of (\ref{DE}-\ref{Ads}) is derived by Hansen \cite{H1961}:
    \begin{equation}    
c_s(t) = c_e-\frac{1}{\sqrt{\pi D}}\int_0^t \frac{1}{(t-\tau)^{1/2}} \frac{d\G(\tau)}{d\tau}\, d\tau.
\label{WTD}
      \end{equation}
This integro-differential relation could be more useful than the Ward-Tordai integral equation in the case of barrier-diffusion control regime, as discussed in \cite{MAF2017}. Indeed, plugging representation (\ref{WTD}) for $c_s(t)$ into the barrier-diffusion control relation (\ref{BDCA}) 
yields a single integro-differential equation for the adsorption $\G(t)$:
\begin{equation}
 \frac{d\G(t)}{dt}=K_a \Phi(\G(t))\left[c_e - \frac{1}{\sqrt{\pi D}}\int_0^t \frac{1}{(t-\tau)^{1/2}} \frac{d\G(\tau)}{d\tau}\, d\tau- f(\G(t))/K \right].
\label{BAFD}
\end{equation}

The approach described in the present work is based on an interpretation of the integrals in the equations (\ref{WTI}-\ref{BAFD}) 
as fractional Riemann-Liouville operators and application of the theory of fractional calculus. 
In eqs. (\ref{WTI}, \ref{WTFI}) they are fractional integrals of order $1/2$ as discussed in \cite{H2016}. 
In eqs. (\ref{WTD}, \ref{BAFD}) the integro-differential operators represent fractional derivatives of order $1/2$. 
The definitions of the fractional operators, as well as some of their basic properties that are used further in the study, are given in the following subsection.

\subsection{Fractional calculus operators}
\label{FCO}

Basic definitions and properties of fractional calculus operators are given next (for more details see e.g. \cite{GM1997}).
%Different types of fractional derivatives and integrals exist, depending on the kernels that they involve, 
%more information can be found in \cite{GM1997}. The most popular, regarding applications to models of different processes, 
%are the Riemann-Liouville fractional operators of order $\al \in (0,1]$ as given below.

The fractional Riemann-Liouville integral $\J_t^\al$ of order $\al>0$ is defined as 
\begin{equation}\label{CFI}
\J_t^\al u(t)=\frac{1}{G(\al)}\int_0^t{(t-\tau)^{\al-1}} {u(\tau)}\,d\tau.
\end{equation}
Here and in what follows $G(\cdot)$ denotes the Gamma function.

The fractional derivative $\D_t^\al$ of order $\al \in (0,1]$ is defined as
\begin{equation}\label{CFD}
\D_t^\al u(t)=\J_t^{1-\al} \frac{d u(t)}{dt}=\frac{1}{G(1-\al)}\int_0^t{(t-\tau)^{-\al}} \frac{d u(\tau)}{d\tau} \,d\tau.
\end{equation}
This definition gives the Caputo fractional derivative, which under the assumption $u(0)=0$ coincides with the Riemann-Liouville fractional derivative. 
The notation $\D_t^1$ for the classical first order derivative is also used in this work.

An important characteristics of the Riemann-Liouville time-fractional operators are their dimensions. 
It can be easily seen that the dimension of the fractional differential and integral operators of 
order $\al$ are $[time^{-\al}]$ and $[time^{\al}]$ respectively:
\begin{equation}\label{FOTrans}
\D_t ^\al u(t)=T_c^{-\al}\D_{\ttt}^\al u(\ttt); \ \ \ \J_t^\al u(t)=T_c^{\al}\J_{\ttt}^\al u(\ttt),\ \ \ t=T_c . \ttt .
\end{equation}

The following relations that hold under the assumption $u(0)=0$ are used in the present study (see e.g. \cite{GM1997}):
\begin{eqnarray}
%&&\D_t^\al u(t) =\D_t^1\J_t^{1-\al}u(t)=\J_t^{1-\al} \D_t^1 u(t),\label{CFDI}\\
&&\D_t^\al \J_t^\al u(t)=\J_t^\al \D_t^\al u(t)=u(t), \ \ \ \ \al\in (0,1)\label{FDFI}\\
&&\D_t^{\al_1} \J_t^{\al_2} u(t)=\D_t^{\al_1-\al_2} u(t),\ \ \ \ 1>\al_1>\al_2>0,\label{FDI}\\ 
&&\J_t^{\al_1} \D_t^{\al_2} u(t)=\J_t^{\al_1-\al_2} u(t),\ \ \ \ 1>\al_1>\al_2>0.\label{DFI}
\end{eqnarray}
 Let us remind that in the present study the case $\G(0)=0$ is considered (which implies also that $f(\G(0))=0$). 
The approach outlined in the present work can be applied also to the case of a nonzero initial adsorption, $\G(0)>0$. 
However, it involves additional terms and parameters in the equations and because of that is not considered here.

In the considerations in the next section the following two-term composite fractional derivative operator appears
\begin{equation}\label{TTFD}
\D_{A,t}^{1,1/2} = A.\D_t^{1} + \D_t^{1/2},
\end{equation}
where $A \ge 0$ is a constant. An integral operator $\J_{A,t}^{1,{1/2}}$ satisfying $$\D_{A,t}^{1,1/2} \J_{A,t}^{1,{1/2}}u(t)=\J_{A,t}^{1,{1/2}}\D_{A,t}^{1,1/2}u(t)=u(t),$$ provided $u(0)=0$, is defined as (see e.g. \cite{BB2022}):
\begin{equation}\label{TTFI}
\J_{A,t}^{1,{1/2}} u(t)=\int_0^t \xi_A(t-\tau) {u(\tau)}\,d\tau
\end{equation}
with kernel 
\begin{equation}\label{TTK}
\xi_A(t) = \frac{1}{A} E_{1/2} \left(-\frac{t^{1/2}}{A} \right) = \frac{1}{A} \exp\left(\frac{t}{A^2}\right) \er \left(\frac{t^{1/2}}{A}\right), 
\end{equation}
where $E_{1/2}$ is the Mittag-Leffler function \cite{GM1997}
\begin{equation}\label{MLF}
E_{1/2}(-z)=\exp({z^2}) \er ({z})=\sum_{n=0}^\infty{\frac{(-1)^n z^{n}}{G(n/2+1)}}.
\end{equation}
%and $\er$ is the complementary error function. 
Taking into account the asymptotic behavior of the Mittag-Leffler function (see e.g. \cite{H1961, GM1997}), the following
asymptotic expansions for the kernel $\xi_A(t)$ are deduced
\begin{eqnarray}
&&\xi_A(t)=\frac{1}{A}-\frac{2t^{1/2}}{A^2 \sqrt{\pi}}+O\left(\frac{t}{A^2}\right),\ \ \ \ \ \ \ \frac{t}{A^2} << 1,\label{as1}\\
&&\xi_A(t)=  \frac{1}{\sqrt{\pi}t^{1/2}}  - \frac{A^2}{2\sqrt{\pi}t^{3/2}}+O\left(\frac{A^4}{t^2}\right), \ \ \ \frac{t}{A^2} >> 1. \label{as2}
\end{eqnarray}

It is known from the theory of fractional calculus that ordinary linear fractional equations 
possess analytical solutions and they are given in terms of Mittag-Leffler functions. 
Such linear fractional equations appear in the following section in the case of Henry adsorption isotherm, 
as well as in regard to asymptotic solutions, derived in Section~\ref{SecAs}. 

It is worth noting that all relations for the fractional operators given above and the analytical results obtained in the present study can be derived applying the technique of Laplace transform. However, the use of fractional calculus operators here makes the steps more visible and manageable, and shortens them. 

In the following sections the mathematical model of the adsorption $\G(t)$ is given in terms of fractional derivative operators, 
Section~\ref{SecFDM}, and fractional integrals, Section~\ref{SecFIM}. 

\section{Fractional derivative model}
\label{SecFDM}
The mathematical model for the adsorption as an ordinary fractional differential equation is considered in the present section. 
First, the limiting case of diffusion-control is considered, where the model includes fractional derivative of order $1/2$, see also \cite{H2016}. 
Results in the other limit, of barrier-controlled adsorption, where the process is governed by first order differential equation, 
are also presented for completeness. 
The main results, concerning the case of mixed barrier-diffusion controlled adsorption, are presented in the second subsection. 
In this case the mathematical model includes two differential operators of orders $1$ and $1/2$, respectively.

\subsection{Limiting cases of diffusion-controlled and barrier-controlled adsorption}
\label{FDMDC}

In the diffusion-controlled regime the adsorption is governed by the Ward-Tordai integral equation (\ref{WTFI}),  
which can be written as a fractional integral equation of order $\al=1/2$ 
    \begin{equation}    
\G(t) = \sqrt{D}\J^{1/2}_t \left[c_e-f(\G(t))/K\right],
\label{DCAIE}
      \end{equation}
taking into account the definition (\ref{CFI}) and the identity $G(1/2)=\sqrt{\pi}$.
Now, applying fractional differentiation $\D^{1/2}_t$ and using (\ref{FDFI}), 
(\ref{DCAIE}) is transformed to a fractional differential equation of order $1/2$:
    \begin{equation}    
\D^{1/2}_t \G(t) = \sqrt{D}(c_e-f(\G(t))/K),\ \ \ \G(0)=0.
\label{DCADE}
      \end{equation}
It is in fact Hansen equation ({\ref{WTD}), with $c_s(t)$ given by (\ref{DCA}), written in terms of fractional derivative (\ref{CFD}). 
Alternatively, applying fractional integration of order $1/2$ to eq.(\ref{WTD}) Ward-Tordai integral equation (\ref{WTI}) is obtained. 
This illustrates how fractional calculus operators appear naturally in the simplest model of adsorption.

The main advantage of the compact representations (\ref{DCAIE}) and (\ref{DCADE}) of the mathematical model (\ref{DE}-\ref{Ads}) and (\ref{DCA}) 
is the reduction of the dimension, the number of unknowns and equations respectively to a single ordinary equation for $\G(t)$. 
Another important advantage is that the number of parameters ($c_e$, $K$, $K_a$, $\G_\infty$, $D$, $\beta$) of the model can be reduced 
by applying a proper transformation to $\G$ and $t$. 
Different transformations are used in the literature, see for instance \cite{DVK2002, KS2023}. 
Here we adopt the characteristic times used in \cite{DVK2002} (see also \cite{DHCL2000}).
The diffusion time in the case of Henry adsorption isotherm is used for nondimensionalization of time:
    \begin{equation}    
t= T_d .\ttt, \ \ \ T_d=\frac{(K\Gi)^2}{D}.
\label{TransT}
      \end{equation}
 The adsorption $\G$ is scaled with its equilibrium value $\G_e=\G(t \to \infty)$:
    \begin{equation}    
\G(t)= \G_e . \G^*(\ttt).
\label{TransG}
      \end{equation}
Using the above scaling, the values of $\G^*$ are in the interval $[0,\, 1]$ ($\G^*(0)=0;\, \G(\infty)=1$), 
which is useful for the graphical presentation of the results. 
For the sake of brevity of the notations a dimensionless parameter $\tb$ is introduced 
    \begin{equation}    
\tb= \beta \frac{2\Gi}{kT}.
\label{TransB}
      \end{equation}
The equilibrium adsorption $\G_e$ can be defined as the solution in the interval $0<\G_e<\Gi$ of the equation:
    \begin{equation}    
f_e = f(\G_e)=K . c_e.
\label{Ge}
      \end{equation}
The above equation has exactly one solution if $\tb < 4$ in the case of Frumkin and $\tb < 6.75$ for van der Waals models. 
Above these values of $\tb$ the corresponding models predict negative surface elasticity.

Transformations different from (\ref{TransT}-\ref{TransG})  are also used in the present work. 
Justification of the choice of transformations will be discussed in the course of the presentation of the study. 
The present choice ($\Gamma^*, \ttt$) has the advantage to reduce the number of parameters in the case of Henry model to a minimum 
($0$ for the diffusion-controlled regime and $1$ in the case of mixed barrier-diffusion controlled adsorption). 

%{\it Remark.} 
Let us note that in this work the notation $\cdot^*$ ($\G^*,\, t^*$) is used when the corresponding characteristic value is defined based on $\G_e$.
If it is based on $\Gi$ the corresponding dimensionless variables are denoted by $\tilde \cdot$ ($\tG, \, \ttt$).

The transformation (\ref{TransT}-\ref{TransG}) is applied to the fractional derivative variant of the Ward-Tordai equation. 
Taking into account the dimensions of the fractional order operators (\ref{FOTrans}) the fractional differential equation (\ref{DCADE}) in dimensionless form reads: 
    \begin{equation}    
\D^{1/2}_{\ttt} \Gs(\ttt) = \Gis \left[f_e - f\left(\Gs\left(\ttt\right).\G_e\right) \right],\ \ \ \Gs(0)=0,
\label{DCADS}
      \end{equation}
where $\Gis$ is the dimensionless maximum adsorption $\Gis=\Gi/\G_e$. 
As it was mentioned, one of the advantage of this transformation is the reduction of the number of parameters. 
Thus, the above equation depends on one parameter $\Gis$. 
Note that for a given adsorption model ($f(\G)$) the values of $f_e=f(\G_e)$ and $f(\G^*(\ttt).\G_e)$ are defined by $\Gis$ 
and vice versa, at given $f_e$, $\Gis$ is determined via eq. (\ref{Ge}). 
In the cases of Frumkin and van der Waals models additional dimensionless parameter $\tb=\Gi\, 2\beta/(k.T)$ appears.

In the case of Henry adsorption model $f(\G(t))=\G(t)/\Gi$ (see Table~\ref{Table1}) and therefore $\Gis\, f(\G_e) =1$ and $\Gis\, f(\Gs(\ttt). \G_e)= \Gs(\ttt)$. Thus, the dimensionless fractional derivative model of diffusion-controlled adsorption (\ref{DCADS}) in the case of Henry isotherm becomes parameterless: 
    \begin{equation}    
\D^{1/2}_{\ttt} \Gs(\ttt) = 1 - \Gs(\ttt),\ \ \ \Gs(0)=0.
\label{DCADHS}
      \end{equation}
For the sake of brevity the zero initial condition, $\Gs (0)=0$, will be omitted further.
			
From the fractional calculus theory there exists a unique solution of equation (\ref{DCADHS}), given in terms of Mittag-Leffler function (\ref{MLF}). 
Applying different representations of this function (see e.g. \cite{GM1997}), the solution can be written in three different forms
\begin{eqnarray}
  	&& \Gs_h(\ttt) = 1- E_{1/2}\left(-\sqrt{\ttt}\right) = 1 -  \exp(\ttt) \er\left(\sqrt{\ttt}\right), \label{SDCAH2} \\ 
	&& \Gs_h(\ttt) = \sum_{n=0}^\infty { \frac{(-1)^n .\ttt^{(n+1)/2}}{G ( (n+3)/2) } }, \ \ \ G(.) - \mbox{Gamma\ function},\label{SDCAH3}\\
	&& \Gs_h(\ttt) = 1 - \int_0^\infty { \frac{\exp(-r.\ttt)}{\pi(r+1)\sqrt{r}} }\,dr= 
	1 - 2\int_0^\infty {\frac{\exp(-r^2.\ttt)}{\pi(r^2+1)} }\, dr. \label{SDCAH1} 
	\end{eqnarray}
%where $G(.)$ is the Gamma function.	
	
The first two forms of the solution $\Gs_h(\ttt)$ are well known (see e.g. \cite{S1952,MB2006,H2016}). 
The series representation (\ref{SDCAH3}) is suitable for asymptotic analysis at small $\ttt$. Both (\ref{SDCAH2}) and (\ref{SDCAH3}) are not suitable for computing %the values of the adsorption 
$\Gs_h(\ttt)$ at large $\ttt$. 
%Similar is the situation with the second form  because of the presence of $\exp(\ttt)$. 
The integral representations in (\ref{SDCAH1}) have the advantage to give accurate numerical results without limitations on the value of $\ttt$.
In our numerical experiments the infinite integrals are computed using MatLab subroutine for numerical integration. 

The other limiting case of barrier-controlled adsorption is also well investigated (see for instance \cite{A1987,MS2015,KDD2016}).
It is achieved at high initial surfactant concentration $c_e$ or fast diffusion, compared to the adsorption. 
Then the subsurface concentration $c_s(t)$ can be assumed constant (independent of time, $c_s(t)=c_s(0)=c_e$). 
Some well-known results, which are used further in this study,  are presented next.
In the barrier-controlled regime the adsorption $\G(t)$ is independent of the diffusion 
equation in the bulk and is governed by a first order ordinary differential equation (equation (\ref{BDCA}) with $c_s(t)=c_e$):
\begin{equation}
 \frac{d\G(t)}{dt}=K_a \Phi(\G(t))\left[c_e - f(\G(t))/K\right].
\label{BCA}
\end{equation}
This equation is made dimensionless using the scaling of the adsorption given in (\ref{TransG}), while time is scaled by the adsorption time:
    \begin{equation}    
t= t^* .T_a , \ \ \ \ \ T_a=\frac{\G_e}{K_a.c_e}.
\label{TransTa}
      \end{equation}
Thus, equation (\ref{BCA}) in terms of $\Gs(t^*)$ reads:
\begin{equation}
 \frac{d\Gs(t^*)}{d t^*}=\left[1 - f(\Gs(t^*).\G_e)/f_e\right]\Phi(\Gs(t^*).\G_e).
\label{BCAD}
\end{equation}

The cases of Henry and Langmuir kinetic models, where well-known analytical solutions exist (see e.g. \cite{A1987}), are considered below.

For Henry model, taking into account the relations
$$f_e=\frac{1}{\Gis},\ \ f(\Gs.\G_e)=\frac{\Gs}{\Gis},\ \ \Phi(\Gs.\G_e)=1, \ \ \G_e = K.c_e. \Gi,$$
equation (\ref{BCAD}) reduces to the parameterless equation
\begin{equation}
 \frac{d\Gs(t^*)}{d t^*}=1 - \Gs(t^*) 
\label{BCAH}
\end{equation}
with analytical solution
\begin{equation}
 \Gs(t^*)= 1-\exp(-t^*). 
\label{SBCAH}
\end{equation}
The adsorption time defined by (\ref{TransTa}) in the case of Henry kinetic model ($\G_e=K.c_e.\Gi$) is 
$$T_a=\frac{K.\Gi} {K_a}.$$

In the case of Langmuir model, taking into account that 
$$f_e=K.c_e = \frac{1}{\Gis-1}, \ \ f(\Gs.\G_e).\Phi(\Gs.\G_e)=\frac{\Gs}{\Gis}, \ \ \G_e = \frac{K.c_e.\Gi}{1+K.c_e} = \frac{f_e. \Gi}{1+f_e},$$ 
the governing dimensionless equation (\ref{BCAD}) and its solution are the same as in the case of Henry model, eq.(\ref{BCAH}-\ref{SBCAH}), 
with the only difference in the scaling of the time $t$: here the adsorption time is additionally scaled by $1/(1+K.c_e)$:
$$T_a=\frac{K.\Gi} {K_a(1+K.c_e)}.$$

\subsection{Mixed barrier-diffusion control}
\label{FDMBDC}

As it was shown at the end of Section \ref{subsecMM} the mathematical model (\ref{DE}-\ref{BDCA}) describing the barrier-diffusion controlled 
adsorption can be reduced to a single integro-differential equation (\ref{BAFD}) for the adsorption $\G(t)$. 
Taking into account the definition of fractional derivative (\ref{CFD}), equation (\ref{BAFD}) can be written, 
after simple rearrangement, as the following composite fractional differential equation
    \begin{equation}    
\frac{\sqrt{D}}{K_a\Phi(\G(t))}\, \D^1_t \G(t) + \D^{1/2}_t \G(t) = \sqrt{D}\left[c_e-f(\G(t))/K\right].
\label{FBDCA}
      \end{equation}
%In terms of fractional calculus classification (\ref{FBDCA}) is a two-term (of orders $1$ and $1/2$) time-fractional differential equation.
To write it in dimensionless form the transformation (\ref{TransT}-\ref{TransG}), used in the diffusion-controlled case, is applied.
In terms of $\Gs(\ttt)$ equation (\ref{FBDCA}) reads
    \begin{equation}    
\frac{\Ba}{\Phi(\Gs(\ttt)\G_e)}\, \D^1_{\ttt} \Gs(\ttt) + \D^{1/2}_{\ttt} \Gs(\ttt) = \Gis\left[f_e-f(\Gs(\ttt).\G_e)\right],
\label{DBDCA}
      \end{equation}
where the dimensionless group $\Ba=D/(K_a .\Gi .K)$ represents the influence of the "Barrier part" of the adsorption to the total adsorption. 
%In the limiting cases $\Ba=0$ and $\Ba=\infty$ equation (\ref{DBDCA}) reduces to that of diffusion-controlled and barrier-controlled adsorption, respectively. 
The dimensionless group $\Ba$ can be considered as the ratio between the adsorption time and the diffusion time in the case of Henry model ($\Ba=T_a/T_d$).

Comparing the above model of mixed adsorption with that in the diffusion-controlled regime (\ref{DCADS}), an 
additional term (with additional parameter $\Ba$) appears here, that corresponds to the barrier part of the adsorption. 
Thus, the parameters of the dimensionless model (\ref{DBDCA}) for the considered 
adsorption isotherms are $\Ba$ and $f_e$, with third parameter $\tb$ in the cases of Frumkin and van der Waals models.

Depending on the function $\Phi(\G)$ (see Table \ref{Table1}) the considered kinetic models are divided in two groups: 
Henry, Volmer and van der Waals, where $\Phi(\G)=1$ and Langmuir and Frumkin with $\Phi(\G)=1-\G/\Gi$. 
From mathematical point of view the main difference between the two groups is with respect to the linearity of the governing equations. 
For the first group, $\Phi(\G)=1$, the operator at the left-hand side of (\ref{DBDCA}) is a linear differential operator, 
while in the second group it is nonlinear. 

To the end of this section the case $\Phi(\G)=1$ is considered, in which equation (\ref{DBDCA}) reads:
    \begin{equation}    
\Ba\, \D^1_{\ttt} \Gs(\ttt) + \D^{1/2}_{\ttt} \Gstt =%\D^{1,1/2}_{\Ba,\ttt}\, \Gstt = 
\Gis\left[f_e-f\left(\Gstt.\G_e\right)\right].
%,\ \ \ \G^*(0)=0
\label{DBDCAF1}
      \end{equation}
The main attention is payed to the case of Henry isotherm where (\ref{DBDCAF1}) reduces to the equation 
    \begin{equation}    
\Ba\, \D^1_{\ttt} \Gstt + \D^{1/2}_{\ttt} \Gstt = 1-\Gstt,
%,\ \ \ \G^*(0)=0
\label{HBDCA}
      \end{equation}
 depending on only one parameter - the dimensionless group $\Ba$. 
Equation (\ref{HBDCA}) possesses analytical solution given in terms of Mittag-Leffler functions, see \cite{GM1997}. 
As in the case of diffusion-controlled adsorption, the solution of equation (\ref{HBDCA}) is given in different forms:

- Via Mittag-Leffler function $E_{1/2}$ (see also Sutherland \cite{S1952}):
\begin{eqnarray}
	&&\Gs_H(\ttt) = 1 +  \frac{1}{\Ba\, (\lambda_+ - \lambda_-)}\left[\frac{1}{\lambda_+} E_{1/2}\left(-\lambda_+\sqrt{\ttt}\right) 
	- \frac{1}{\lambda_-} E_{1/2}\left(-\lambda_-\sqrt{\ttt}\right) \right] \label{SBDCAH2}\\ 
	&&= 1 +  \frac{1}{\Ba\, (\lambda_+ - \lambda_-)}\left[\frac{1}{\lambda_+} \exp\left(\lambda_+^2 \ttt\right) \er\left(\lambda_+ \sqrt{\ttt}\right) 
	- \frac{1}{\lambda_-} \exp\left(\lambda_-^2 \ttt\right) \er\left(\lambda_- \sqrt{\ttt}\right) \right], \nonumber
      \end{eqnarray}
			where $\lambda_{\pm}$ are the roots of the equation %corresponding to (\ref{HBDCA}) quadratic equation in the Laplace domain, 
			$\Ba \, x^2 + x + 1 = 0$:
			$$\lambda_{\pm} = -\frac{1}{2\Ba} \pm \sqrt{\frac{1}{4\Ba^2}-\frac{1}{\Ba}}.$$
			
- As infinite sums: 
\begin{eqnarray}
	&&\Gs_H(\ttt) = \frac{1}{\Ba} \sum_{n=1}^\infty{ \frac{(\lambda_+^n- \lambda_-^n)}{(\lambda_+ - \lambda_-)}\frac{\ttt^{\frac{n+1}{2}}}{G \left( \frac{n+3}{2}\right) }}, \label{SBDCAH3}\\
	&&\Gs_H(\ttt)=\sum_{n=0}^\infty \sum_{m=0}^\infty{ \frac{(n+m)!}{n!m!}\frac{(-1)^{n+m}}{\Ba^{n+m+1}} \frac{\ttt^{n+m/2+1}}{G \left(n+m/2+2\right)} },\label{MML}
	%&&\Gs_H(\ttt)=\sum_{n=0}^\infty \sum_{m=0}^\infty{ \frac{(n+m)!}{n!m!}\frac{(-1)^{n+m}}{\Ba^{n+m+1}} \frac{\ttt^{\frac{2n+m+2}{2}}}{G (\frac{2n+m+4}{2})} }.\nonumber
      \end{eqnarray}
			where $G(.)$ is the Gamma function.	
			
			- Integral representations (see also \cite{BB2022}):
    \begin{equation}    
	\Gs_H(\ttt) = 1 - \int_0^\infty \frac{\exp(-r.\ttt)\, dr}{\pi\sqrt{r}((\Ba .r-1)^2 + r)} 
	= 1 - \frac{2}{\pi} \int_0^\infty {\frac{\exp(-r^2.\ttt)\, dr}{(\Ba .r^2-1)^2 + r^2}}. 
	\label{SBDCAH1}
      \end{equation}
			
The series expansion (\ref{MML}) is derived from the representation of the solution of equation (\ref{HBDCA}) in terms of multinomial Mittag-Leffler function 
(see e.g. \cite{BB2022,FCAA2021}).

In analogy to the solution $\Gs_h(\ttt)$ (see (\ref{SDCAH2}-\ref{SDCAH1})) in the case of diffusion-controlled adsorption, 
the above forms of solution $\Gs_H(\ttt)$ in the mixed adsorption regime can be useful in different situations. 
The integral representations (\ref{SBDCAH1}) have the advantage to give numerical results for all values of $\ttt > 0$. 
Another advantage of (\ref{SBDCAH1}) is that it directly reduces to $\Gs_h(\ttt)$ given by (\ref{SDCAH1}) when $\Ba=0$, 
while the other two forms (\ref{SBDCAH2}-\ref{SBDCAH3}) are not defined for $\Ba=0$. 

\begin{figure}[h]%% placement specifier
\centering%% For centre alignment of image.
\includegraphics[width=12cm]{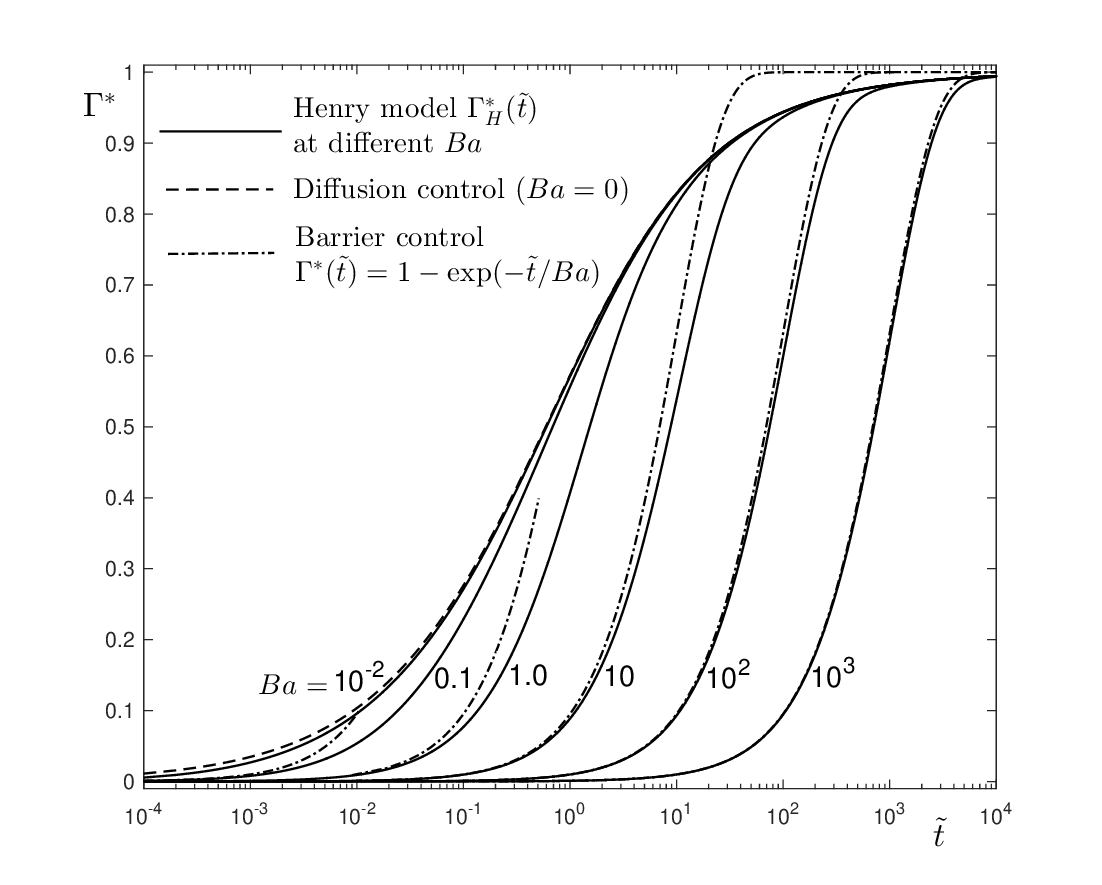}
\caption{Barrier-diffusion controlled adsorption $\Gs_H(\ttt)$ (eq.(\ref{SBDCAH1}) - solid lines) for Henry kinetic model at different $\Ba$. The limiting cases of diffusion control and barrier control are given by dashed and dash-dotted lines respectively.}
\label{Fig1}
\end{figure}
In Fig.~\ref{Fig1} adsorption is given at different values of $\Ba$ in a broad range $10^{-2}\le \Ba \le 10^3$. The results presented by the solid lines are
based on (\ref{SBDCAH1}) computed using subroutine for numerical integration in MatLab. 
The results indicate that at $\Ba \le\ 10^{-2}$ the adsorption is already diffusion-controlled, while for $\Ba \ge\ 10^2$ it is barrier-controlled. 
It is also seen that at moderate $\Ba$ and for small surface coverage ($\Gs<<1$) the process is barrier-controlled, while at the end ($\Gs$ close to $1$) it is diffusion-controlled.

For the other representations % (\ref{SBDCAH2}-\ref{SBDCAH3}) 
there are computational problems at high values of $\ttt/\Ba$ ($\Gs$ close to $1$). 
Also, the form (\ref{SBDCAH2}) involves complex values of the roots $\lambda_\pm$ at $\Ba>1/4$, 
which introduces additional challenges for the numerical computation. 
The series representations are suitable for analysis of the asymptotic behavior at $\ttt/\Ba<<1$. 
The first terms of (\ref{MML}) are given below:
\begin{eqnarray} 
	\Gs_H(\ttt) &=& \frac{1}{\Ba} \left[ \ttt + \frac{4}{3}\frac{\ttt^{3/2}}{\sqrt{\pi}\, \Ba} + \frac{\ttt^2}{2} \left(\frac{1}{\Ba^2} - \frac{1}{\Ba}\right) +
	\frac{8 \, \ttt^{5/2}}{15\, \sqrt{\pi}} \left(\frac{2}{\Ba^2} - \frac{1}{\Ba^3}\right) \right] \nonumber  \\
	&+& O\left(\left(\frac{\ttt}{\Ba}\right)^3\right).\label{AsH}	
      \end{eqnarray}
It should be noted that in the case of adsorption relaxation \cite{DVK2002} the authors have derived very similar expressions 
to the integral representation (\ref{SBDCAH1}) and the above asymptote (see their eqs.(13-15)). 
The difference between the expressions here and in \cite{DVK2002} is due to differences in the scaling and the initial condition.

An important characteristics of the solution $\Gs_H(\ttt)$ of the Henry model 
is that it depends on material parameters only ($\Ba=D/(K_a.\G_\infty.K);\ \ttt=t.D/(K.\G_\infty)^2$). 
An advantage is also that it depends only on one parameter - $\Ba$.
These properties of $\Gs_H(\ttt)$ are essential in regard to the second order (at small surface coverage) 
asymptotes for the adsorption and the surface tension derived in Section~\ref{SecAs}.

\section{Fractional integral model}
\label{SecFIM}

The transformation of fractional differential to fractional integral model is performed by applying an 
integral operator, which is inverse to the corresponding fractional differential one in the equation for the adsorption. 
For instance, the Ward-Tordai integral equation in dimensionless form is obtained from the fractional derivative model (\ref{DCADS}) 
of diffusion-controlled adsorption by applying fractional integral $\J^{1/2}_{\ttt}$, which yields by the use of  (\ref{FDFI})
    \begin{equation}    
\Gs(\ttt) = \Gis \J^{1/2}_{\ttt}\left[f_e - f(\Gs(\ttt).\G_e) \right]= \int_0^{\ttt} \frac{\Gis[f_e-f(\Gs(\tau).\G_e)]}{\sqrt{\pi}(\ttt-\tau)^{1/2}}\, d\tau.
\label{DCAIS}
      \end{equation}
			
This approach is used in the following subsections where the two variants of the function $\Phi(\G)$ %that appears in the considered kinetic models 
(see Table~\ref{Table1}) are analyzed separately.

\subsection{The case $\Phi(\G)=1$}
\label{SubSecFIM1}

In the case of Henry, Volmer or van der Waals adsorption isotherms ($\Phi(\G)=1$) the left-hand side of 
the fractional derivative model (\ref{DBDCA}) is linear, see (\ref{DBDCAF1}). 
For this case two possibilities exist to transform the fractional derivative equation to fractional integral one. 

The first approach is to consider the composite fractional derivative operator $\D^{1,1/2}_{\Ba,\ttt}$, defined by (\ref{TTFD}), and %In the second option the two derivatives of order $1$ and $1/2$ are treated separately. 
 apply its inverse integral operator $\J^{1,1/2}_{\Ba,\ttt}$, defined by (\ref{TTFI}). In this way the equation (\ref{DBDCAF1}) is transformed into the following integral equation 
\begin{eqnarray}
\Gs(\ttt) &= &\J^{1,1/2}_{\Ba,\ttt} \left\{\Gi^*\left[f_e-f(\G^*(\ttt).\G_e)\right]\right\} \nonumber \\ 
&=&\int_0^{\ttt}\xi_\Ba (\ttt-\tau)\, \Gis\left[f_e-f(\G^*(\tau).\G_e)\right]\, d\tau, \label{IDCAF1}  
%&=\int_0^{\ttt}{\Ba^{-1}\exp \left((\ttt-\tau)/\Ba^2\right)\, \er\left(\sqrt{\ttt-\tau}/\Ba\right) \, \Gi^*\left[f_e-f \left(\Gs(\tau).\G_e \right)\right]}\, d\tau \nonumber
%,\ \ \ \G^*(0)=0
      \end{eqnarray}
			where the kernel $\xi_\Ba$ is defined in (\ref{TTK})
			$$
			\xi_\Ba(t)=\Ba^{-1}\exp \left(t/\Ba^2\right)\, \er\left(\sqrt{t}/\Ba\right). 
			$$
%Taking into account that $\lim_{x\to \infty} \er(x)  = \exp(-x^2)/(\sqrt{\pi}.x)$ (see for instance \cite{erfc}) it follows that:
    Asymptotic expansion (\ref{as2}) implies
		\begin{equation}    
\left. \xi_\Ba (t) \right|_{Ba=0} = \frac{1}{\sqrt{\pi t}},
\label{KernBa}
      \end{equation}
which means that at $\Ba=0$ (\ref{IDCAF1}) reduces to the Ward-Tordai equation (\ref{DCAIS}). 
Therefore, the integral equation (\ref{IDCAF1}) is a generalization of the well-known Ward-Tordai equation to the case of mixed barrier-diffusion controlled adsorption.	
%This integral equation can be considered as generalization of the Ward-Tordai integral equation to the case of mixed barrier-diffusion controlled adsorption. 

In the next section it will be shown that the well-known first order asymptote for the barrier-diffusion controlled adsorption 
can be directly obtained from (\ref{IDCAF1}). 
Also based on (\ref{IDCAF1}), the existing techniques for numerical integration of the Ward-Tordai equation (see for instance \cite{CCE2010,N2022})
can be applied directly to the case of mixed barrier-diffusion controlled adsorption. 

In the second approach we apply the first order integral operator $\J^1_{\ttt}$ to the differential equation (\ref{DBDCAF1}), which yields
    \begin{equation}    
\Gstt = \frac{1}{\Ba}\left[ \Gis \J^{1}_{\ttt}\left(f_e - f(\Gstt.\G_e)\right) - \J^{1/2}_{\ttt}(\Gstt)\right],
\label{IBDCAF1}
      \end{equation}
where it is taken into account that $\J^1_{\ttt}\D^{1/2}_{\ttt}=\J^{1/2}_{\ttt}$, see (\ref{DFI}). 

Let us denote the right hand side of the above equation by $\F(\Gs)$:
    \begin{equation}    
\F(\Gstt) = \frac{1}{\Ba}\left[ \Gis \int_0^{\ttt}{\left(f_e - f(\Gs(\tau).\G_e)\right)\, d\tau} - \frac{1}{\sqrt{\pi}}\int_0^{\ttt}{\frac{\Gs(\tau)}{\sqrt{\ttt-\tau}}\, d\tau}\right].
\label{DFG}
      \end{equation}
This notation will be convenient for the following case.

\subsection{The case $\Phi(\G)=1-\G/\Gi$}
\label{SubSecFIMF}

The approach in the cases of Langmuir and Frumkin adsorption isotherms is similar. 
However, due to the non-linearity of the first term in (\ref{DBDCA}), additional manipulations (similar to these used in \cite{LRP1996}) are necessary.

Applying integral operator $\J^1_{\ttt}$ to the fractional derivative model (\ref{DBDCA}) and taking into account that
   $$% \begin{equation}    
\int_0^{\ttt}{\left[\frac{1}{\Phi(\Gs(\tau)\G_e)}\frac{ d \Gs(\tau)}{d\tau}\right]\, d\tau} = 
\int_0^{\ttt}{\frac{ d \Gs(\tau)}{1-\Gs(\tau)/\Gis}} = -\Gis \ln(1- \Gstt/\Gis) 
%= F(\G^*) 
%\label{IFG}
 $$  %   \end{equation}
under the assumption in the present study $\G^*(0)=0$, we deduce the
fractional integral equation for the adsorption $\Gs(\ttt)$ in the cases of Langmuir and Frumkin adsorption models 
    \begin{equation}    
\Gstt = \Gis \left\{ 1-\exp\left[-\F(\Gstt)/\Gis\right] \right\},
\label{IBDCAFG}
      \end{equation}
where $\F(\Gs)$ is defined in (\ref{DFG}).

The integral equation (\ref{IBDCAFG}) appears to be very useful for deriving asymptotes of higher order in the case of 
barrier-diffusion controlled adsorption, which is discussed in the following section.
%in the case of Langmuir and Frumkin kinetic models, 

In the presented in this section integral models (\ref{IBDCAF1}) and (\ref{IBDCAFG}) the adsorption evolution is governed by Volterra type integral equations. 
They are a good starting point to develop numerical technique for computer simulation of the adsorption process. 
It is worth noting that the level of difficulty regarding numerical approximation of the integrals in these equations is determined by 
the integral of order $1/2$ which is singular of convolution type. 
Thus, the numerical difficulties for solving the fractional integral models 
of mixed barrier-diffusion controlled adsorption are not grater than these for the classical Ward-Tordai integral equation. 
  
\section{Higher order asymptotes for the adsorption and surface tension}
\label{SecAs}

A significant part of theoretical and experimental studies of the evolution of the adsorption and respectively the surface tension 
use asymptotic predictions (see for instance \cite{A1987, MB2006, BMDSC2024}).
The most popular, so called short-time asymptotes, are based on the assumption of desorption free process. 
In spite of the name "short-time" these asymptotes are used for time of order up to $10^4$ $[s]$, see for instance 
\cite{LHSL2022} in the case of diffusion-controlled and \cite{LPL2022} in the case of barrier-diffusion controlled adsorption. 
In the present section the asymptotes are analyzed in connection with approximations of the desorption term, function $f(\G)$.

\subsection{Approximation of the desorption term}
\label{SSecAsDes}

As a first step of the derivation of higher order asymptotes of the adsorption we consider approximations of the function $f(\G)$ 
that is part of the adsorption isotherm in the considered kinetic models given in Table~\ref{Table1}. 

Taking into account the Taylor expansions of the functions
$$u/(1-u)=u+u^2+u^3+...\ \ \mbox{and}\ \  \exp(u)=1+u+u^2/2+u^3/6+...,$$ 
that are part of the adsorption isotherms,
it is easy to obtain approximation of $f(\G)$ at small values of $\G$ with required order of approximation. 
The leading two terms can be given in a common form for all considered models:
    \begin{equation}    
f(u) = u/\Gi + (Q-\bt)\left(\frac{u}{\Gi}\right)^2 + O\left(\left(\frac{u}{\Gi}\right)^3\right),
\label{TE}
      \end{equation}
where the constant $Q$ depends on the kinematic model under considerations 
\begin{eqnarray}
Q&=&0\mbox{ - Henry,}\nonumber\\
Q&=&1\mbox{ - Langmuir and Frumkin,}\nonumber\\
Q&=&2\mbox{ - Volmer and van der Waals.}\nonumber
\end{eqnarray}
Having the above in mind, and also that $\Gis=\Gi/\G_e$, the right-hand side of the fractional derivative models (\ref{DCADS}) of diffusion-controlled 
and (\ref{DBDCA}-\ref{DBDCAF1}) of barrier-diffusion controlled adsorption reads
\begin{eqnarray}
&R(\Gs(\ttt)) &= \Gis(f_e-f(\Gs(\ttt).\G_e)) \label{RHS}\\
&&=\Gis.f_e - \Gs(\ttt)- \Gis(Q-\bt)\left(\frac{\Gs(\ttt)}{\Gis}\right)^2 + O\left(\left(\frac{\Gs(\ttt)}{\Gis}\right)^3\right).\nonumber
\end{eqnarray}
Note that $R(\Gs(\ttt))$ appears also in the fractional integral models (\ref{DCAIS}-\ref{IDCAF1}). 
At high values of the ratio $\Gis=\Gi/\G_e$ the right-hand side $R(\Gs)$ approaches that of the corresponding Henry model:
$$R_H(\ttt)=\lim_{\Gis \to \infty}R(\Gs(\ttt))=1 - \Gs(\ttt),$$%%\ \ \ (\mbox{because} \ \lim_{\Gis \to \infty}\Gis.f_e =1)$$
where we have used $\lim_{\Gis \to \infty}\Gis.f_e =1.$
This indicates that with decreasing the surfactant concentration ($\G_e \to 0,\ \Gis=\Gi/\G_e \to \infty$) 
%independently on the model (among that considered here, see Table \ref{Table1}) and the parameters in it 
the adsorption $\Gstt$, for all models considered here, tends to that predicted by Henry model $\Gs_H(\ttt)$.

Two approximations of $R(\Gs(\ttt))$ are considered below. First order approximation, with error $O(\Gs)$:
    \begin{equation}    
R_1(\Gs(\ttt)) = \Gis.f_e,
\label{RHS1}
      \end{equation}
which in fact is the well known case of short-time asymptote, where the desorption ($f(\G)$) is neglected.
An estimate of the relative error adopting the above approximation is:
    \begin{equation}    
Er_1 = \frac{| R(\Gs(\ttt)) - R_1(\Gs(\ttt))|}{R_1(\Gs(\ttt))} = \frac{\Gs(\ttt)}{\Gis.f_e} + O\left({\Gs}^2(\ttt)\right).
\label{ERHS1}
      \end{equation}
This shows that $R_1(\Gs)$ is first order approximation of $R(\Gs)$ with respect to $\Gs(\ttt)$, and the leading error is inversely proportional to $\Gis.f_e$.
 
The well known short-time asymptotes of the adsorption (see for instance \cite{MB2006}) can be obtained directly from the generalized Ward-Tordai equation 
(\ref{IDCAF1}). 
Thus, using the first order approximation (\ref{RHS1}) at the right-hand side of (\ref{IDCAF1}), the asymptotic solution in the case of barrier-diffusion controlled adsorption can be obtained. 
In the case of Henry, Volmer and van der Waals models ($\Phi(\Gs)=1$) the asymptote is:
\begin{eqnarray}
&&\Gs(\ttt) \approx \J^{1,1/2}_{\Ba,\ttt}(\Gis .f_e) = \Gis .f_e \int_0^{\ttt} \xi_{\Ba} (\ttt) \, d\tau \nonumber\\ 
&&= \Gis .f_e\left(2\sqrt{\frac{\ttt}{\pi}}+\Ba^2 \xi_{\Ba} (\ttt)-\Ba \right)\label{FOABDC}\\
&&= \Gis .f_e\left(2\sqrt{\frac{\ttt}{\pi}}+\Ba. \exp{\left(\frac{\ttt}{\Ba^2}\right)}.\er{\left({\frac{\sqrt{\ttt}}{\Ba}}\right)}-\Ba\right),\nonumber
\end{eqnarray}
obtained by integration using the definition (\ref{MLF}) of Mittag-Leffler function $E_{1/2}$, respectively the kernel $\xi_{\Ba}$, 
as infinite sum. 
Alternatively, approximation (\ref{FOABDC}) can be obtained using the Laplace transform technique (see the Appendix A in \cite{MB2006}).
In the limit of $\Ba\to 0$ relation 
 (\ref{FOABDC}) yields the well-known short-time asymptote of diffusion-controlled adsorption 
$2. \Gis .f_e \sqrt{\ttt/\pi}$ (see (\ref{KernBa})).
The asymptote (\ref{FOABDC}) is used also in the case of Langmuir model (see for instance \cite{MB2006} and \cite{LPL2022}), 
under the assumption $\G<<\G_\infty$ ($\Phi(\Gs)\approx 1$). 

To take into account the effect of $\Phi(\Gs)=1-\G/\G_\infty$ in Langmuir and Frumkin models, an approach similar to that of \cite{MS2015} is used. 
The operator $\F(\Gstt)$ in (\ref{IBDCAFG}), which is also the right-hand side of (\ref{IBDCAF1}), is approximated using (\ref{IDCAF1}).
This leads to the following asymptote 
\begin{equation}
\Gs(\ttt) \approx \Gis \left\{ 1-\exp\left[f_e\left(2\sqrt{\frac{\ttt}{\pi}}+\Ba. \exp{\left(\frac{\ttt}{\Ba^2}\right)}.\er{\left({\frac{\sqrt{\ttt}}{\Ba}}\right)}-\Ba\right) \right]\right\}. \label{FOABDCG}
\end{equation}

Let us note that such a result is derived using Laplace transform technique in \cite{MS2015} (with a correction \cite{MS2019}). 

The error estimate (\ref{ERHS1}) indicates that the above asymptotes approximate the adsorption with accuracy proportional to $\Gs(\ttt)/(\Gis .f_e)$.
Via the fractional integration in (\ref{IDCAF1}) the error, and respectively the accuracy, depends on time as well. 
However, the first order accuracy, $O(\Gs(\ttt))$, is determining. 
That is why the above asymptotes are applicable at small surface coverage $(\Gs(\ttt)<<1)$, rather than at short time. 
This, together with the dependence of the accuracy on $(\Gis .f_e)^{-1}$, is demonstrated below by comparisons with numerical results.

At the second approximation of $R(\Gs)$ the desorption $f(\Gs(\ttt)\G_e)$ is approximated by its leading term (see (\ref{TE}-\ref{RHS})):
    \begin{equation}    
R_2(\Gs(\ttt)) = \Gis.f_e - \Gs(\ttt).
\label{RHS2}
      \end{equation}
Similar to the first case estimate of the relative error of the above approximations leads to: 
    \begin{equation}    
Er_2 = \frac{| R(\Gs(\ttt)) - R_2(\Gs(\ttt))|}{R_2(\Gs(\ttt))} = \frac{(Q-\bt)({\Gs}^2(\ttt))}{\Gis.f_e.\Gis-\Gs(\ttt)} + O\left({\Gs}^3(\ttt)\right).
\label{ERHS2}
      \end{equation}
This shows that $R_2(\Gs(\ttt))$ is a second order approximation of $R(\Gs(\ttt))$ with respect to $\Gs(\ttt)$ and the leading error of 
the approximation is inversely proportional to $\Gis.f_e.\Gis$.
The above also shows that at $\tb=Q$ the approximation (\ref{RHS2}) of the right-hand side $R(\Gs(\ttt))$ is of third order (with an error $O({\Gs}^3(\ttt)))$.
%An important characteristic of the above approximation is that at fixed value of $\Gis.f_e$ it is universal for all kinetic models considered here. 

\subsection{Second order asymptotes of the adsorption}
\label{SSecAsAds}

Second order asymptotes of the adsorption can be obtained by the approximation (\ref{RHS2}) of the right-hand side of the fractional differential equations 
(\ref{DCADS}) and (\ref{DBDCA}). 
Let us first consider the case of diffusion-controlled adsorption. 
Applying the approximation (\ref{RHS2}) to the model (\ref{DCADS}) leads to linear fractional differential equation of order $1/2$ for the asymptote $\Gs_a(\ttt)$:
    \begin{equation}    
\D^{1/2}_{\ttt} \Gs_a(\ttt)= \Gis .f_e - \Gs_a(\ttt). 
\label{SOADC}
      \end{equation}
The solution of the above equation can be given via the solution $\Gs_h(\ttt)$ of the Henry model (\ref{DCADHS}). 
Indeed, because of the linearity of the fractional operator $\D^{1/2}_{\ttt}$, it is easy to see that the solution of (\ref{SOADC}) is:
    \begin{equation}    
\Gs_a(\ttt)= \Gis .f_e.\Gs_h(\ttt), 
\label{SSOADC}
      \end{equation}
where $\Gs_h(\ttt)$ is given by (\ref{SDCAH2}-\ref{SDCAH1}).

Taking into account that $\Gs_h(\ttt)$ does not include any parameter	the above asymptote depends on one parameter only ($\Gis.f_e$). 
Thus, similar to the first order asymptote considered earlier, the above is applicable for all considered models. 
The effect of a specific model (also parameter $\tb$) is taken onto account via the value of $\Gis.f_e$.	

\begin{figure}[h]%% placement specifier
\centering%% For centre alignment of image.
\includegraphics[width=12cm]{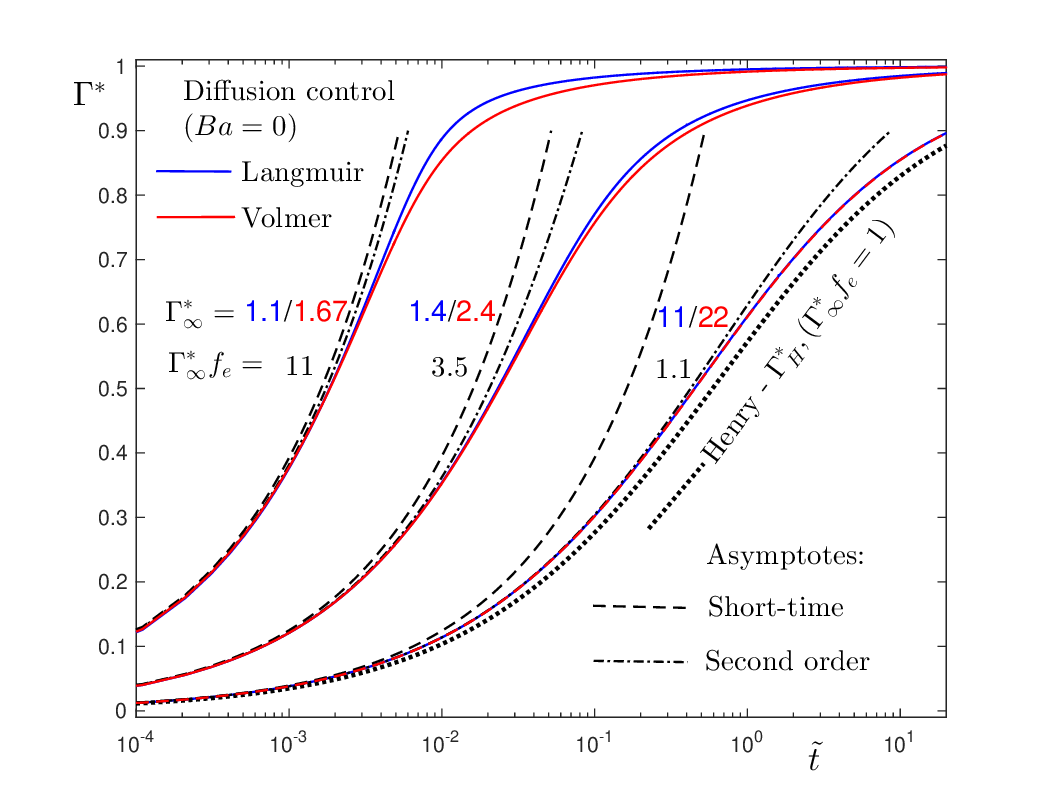}
\caption{Comparison of the short-time (dashed lines) and second order (eq. (\ref{SSOADC}), dash-dotted lines) asymptotes with numerical 
results (solid lines). The presented results are in the case of diffusion-controlled adsorption for Langmuir and Volmer kinetic models at different values of $\Gis .f_e$.}
\label{Fig2}
\end{figure}
In Fig.~\ref{Fig2} short-time and the second order asymptotes (\ref{SSOADC}) are compared with numerical results in the case of 
Langmuir and Volmer kinetic models at different values of $\Gis .f_e$.  
Regarding the comparison between the asymptotes it can be concluded that for relatively high values of $\Gis .f_e$ ($\Gis .f_e=11$) both asymptotes 
give satisfactory prediction in a wide range of the adsorption $\G$, up to about $2/3$ of the equilibrium value $\G_e$. 
With decreasing of $\Gis .f_e$ the accuracy of the short-time asymptote decreases (dashed line), while the accuracy of the second order 
(dash-dotted line) is almost insensitive on the values of $\Gis .f_e$. 
Such difference in the range of applicability of the asymptotes can be explained, except with the difference in the orders of approximation, also with the coefficients of the leading term in the error estimates (\ref{ERHS1}) and (\ref{ERHS2}) respectively. 
For the short-time asymptote this coefficient is $(\Gis .f_e)^{-1}$, while for the second order it is $(\Gis .f_e.\Gis )^{-1}$. 
The numerical results presented in Fig.~\ref{Fig2} indicate that with increasing $\Gis$ ($\Gis .f_e \to 1$) the adsorption $\Gs (\ttt)$ for both, Langmuir and Volmer, models tends to that of Henry model, $\Gs_h$ (dotted line). 

In the general case of mixed barrier-diffusion controlled adsorption the two cases $\Phi(\G)=1$ and $\Phi(\G)=1-\G/\Gi$ 
are considered separately as in the previous section. 

For the simpler case $\Phi(\G)=1$ (Volmer and van der Waals)  the approach described above for the diffusion-controlled regime can 
be applied directly. 
Indeed, approximating the right-hand side of the fractional derivative model (\ref{DBDCAF1}) by (\ref{RHS2}), the equation for the second order asymptote $\Gs_V$ reads:
    \begin{equation}    
\Ba\, \D^{1}_{\ttt} \Gs_V(\ttt) + \D^{1/2}_{\ttt} \Gs_V(\ttt)= \Gis .f_e - \Gs_V(\ttt). 
\label{SOABDC1}
      \end{equation}
This equation is similar to that in the case of Henry model (\ref{HBDCA}). 
The only difference is in the coefficient $\Gis .f_e$ at the right-hand side, which in the case of Henry model is $\Gis .f_e=1$. 
Taking this into account, and the linearity of the equation, its solution $\Gs_V(\ttt)$ is given:
    \begin{equation}    
\Gs_V(\ttt)= \Gis .f_e.\Gs_H(\ttt), 
\label{SSOABDC1}
      \end{equation}
via the solution of the corresponding Henry model $\Gs_H(\ttt)$, see (\ref{SBDCAH2}-\ref{SBDCAH1}). 

\begin{figure}[h]%% placement specifier
\centering%% For centre alignment of image.
\includegraphics[width=12cm]{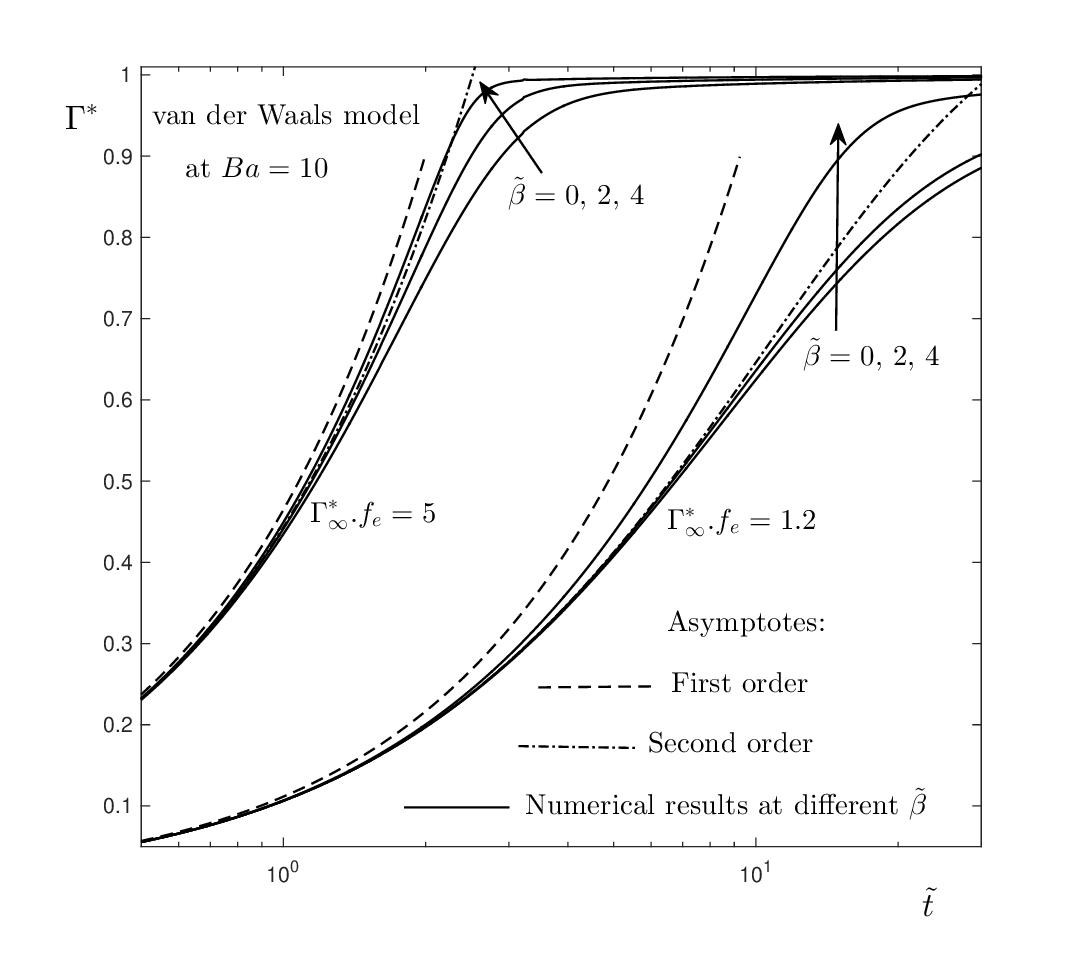}
\caption{Comparison of the numerical results (solid line) with the asymptotes, first (eq. (\ref{FOABDC}), dashed lines) and second order 
(eq. (\ref{SSOABDC1}), dash-dotted lines) for van der Waals model at $\Ba=10$ and different values of $\tb$. 
For the values of $\tb$ (0; 2; 4) the corresponding values of $\Gis$ are: (2.06; 1.54; 1.31) at $\Gis .f_e=5$ and (11.7; 3.35; 1.585) at $\Gis .f_e=1.2$}
\label{Fig3}
\end{figure}

The above second order asymptote of the mixed barrier-diffusion controlled adsorption in the case of van der Waals kinetic model 
$\Gs_V$ is compared with numerical results. 
The comparisons show good agreement in a broad range of the parameters $\Ba$, $\Gis .f_e$ and $\tb$. 
Fig.~\ref{Fig3} presents results at $\Ba=10$, two values of $\Gis .f_e=5;\,1.2$ and different values of the parameter $\tb$.  
The numerical results show that for small surface coverage ($\Gs<<1$) the adsorption $\Gs(\ttt)$ is almost insensitive of the values of $\tb$. 
This is consistent with the fact that the second order approximation of the right-hand side (\ref{RHS2}) depends on $\tb$ only via the value of $\Gis .f_e$. 
At high values of the adsorption, close to the equilibrium $\G_e$ ($\Gs$ close to 1), the influence of the parameter $\tb$ is significant. 
The error estimates (\ref{ERHS2}) shows that in the case of van der Waals model ($Q=2$) the asymptote 
(dash-dotted line) at $\tb=2$ is of third order, which is in agreement with the presented results. 
It should be taken into account that $\G_e$, and also $\Gis=\Gi/\G_e$ depend on $\tb$ (see the caption of the figure). 
The comparisons show that the second order asymptotes (dash-dotted line) give satisfactory results in a broad range of values for the dimensionless 
maximum surfactant adsorption ($1.3<\Gis<12$).
Regarding the comparisons between the first and second order asymptotes, the conclusion from the case of diffusion-controlled adsorption 
(see Fig.~\ref{Fig2}) are confirmed also here.

The case of Langmuir and Frumkin models ($\Phi(\Gs.\G_e)=1-\Gs/\Gis$) is considered below.
It can be expected that (\ref{SSOABDC1}) is a good approximation of $\Gs$ also in this case, especially at $\Gis>>1$ ($\Gi>>\G_e)$ where $\Phi(\G)\approx 1$. 
However, for the solutions of these models $\Gs_V(\ttt)$ is of first order, with error $O(\Gs/\Gis)$.
This is confirmed by comparisons with numerical results at moderate values of $\Gis$, some of which can be seen in Fig.~\ref{Fig4}.
The figure shows that at $\Gis.f_e=11$ ($\Gis\approx 1$) $\Gs_V$ (red dash-dotted line) gives good prediction for Volmer model (solid red line). 
However, the accuracy of its prediction for Langmuir model (solid blue line) is not better than first order asymptote (dashed line). 
That is why it is worth deriving higher order asymptote for the case of Langmuir and Frumkin models. % ($\Phi(\G)=1-\G/\Gi$). 
To do this let us consider the fractional integral models (\ref{IBDCAF1}) and (\ref{IBDCAFG}).
It is seen that for the two cases of $\Phi(\G)$ the models are connected via the fractional operator $\F (\Gs(\ttt))$:
$$\Gs(\ttt)= \F (\Gs(\ttt)) \ \ \mbox{at} \ \ \Phi(\G)=1 \ \ \ \mbox{and}$$
$$\Gs(\ttt) = \Gis \left\{ 1-\exp\left[-\F(\Gs(\ttt))/\Gis\right] \right\} \ \ \mbox{at} \ \ \Phi(\G)=1-\G/\Gi.$$

The asymptote $\Gs_V$ given by (\ref{SSOABDC1}) is of second order for the solution of Volmer and van der Waals models (first equation above, at $\Phi(\G)=1$).
Thus, $\Gs_V$ can be considered as an approximation for $\F (\Gs(\ttt))$, which used in the second equation above leads to:
    \begin{equation}    
\Gs_L(\ttt)= \Gis\left\{1-\exp\left[-f_e.\Gs_H(\ttt)\right]\right\}. 
\label{SSOABDC2}
      \end{equation}

\begin{figure}[h]%% placement specifier
\centering%% For centre alignment of image.
\includegraphics[width=12cm]{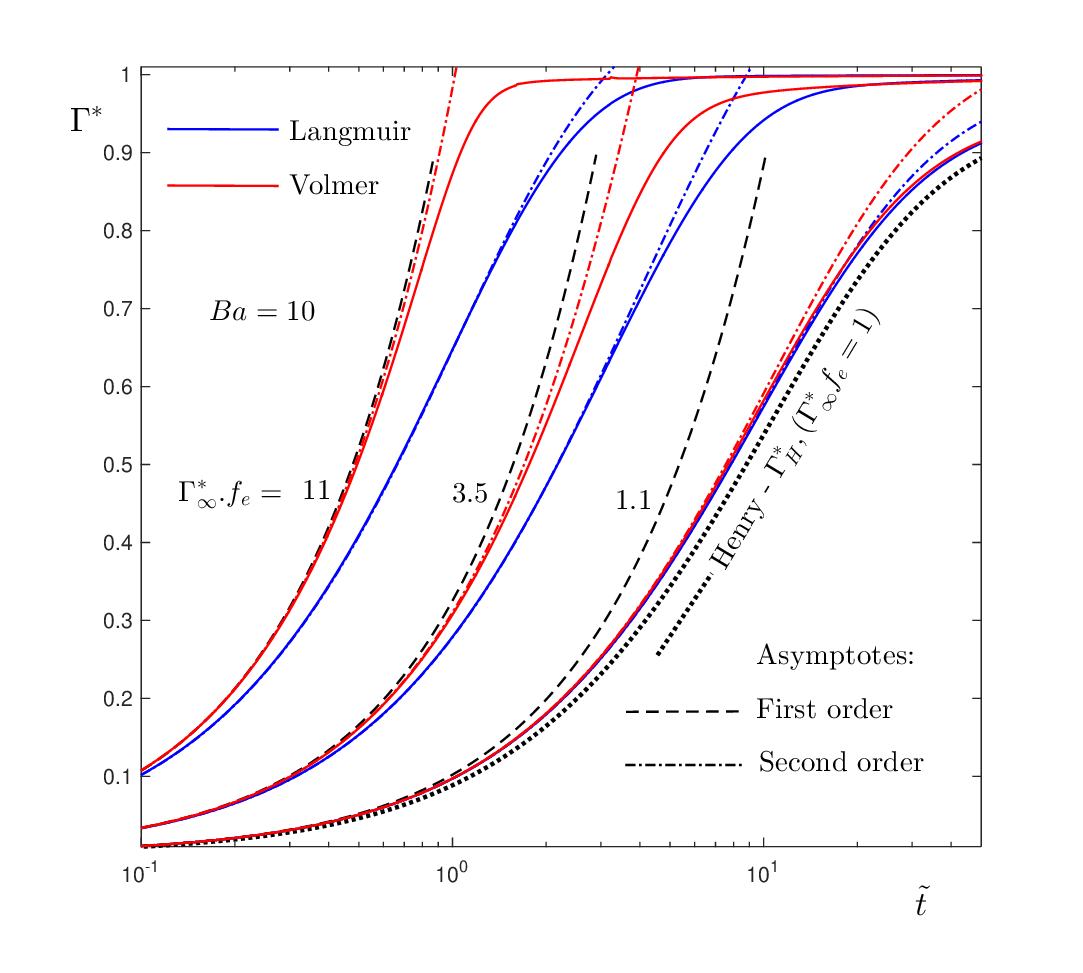}
\caption{Comparison of the second order asymptotes (dash-dotted lines) with numerical results (solid line) at $\Ba=10$ and different $\Gis .f_e$ 
in the case of Langmuir $\Gs_L$ (blue) and Volmer $\Gs_V$ (red). 
The first order asymptotes (eq. (\ref{FOABDC}), dashed lines) as well as the solution $\Gs_H$ of Henry model (eq. (\ref{SBDCAH1}), dotted line) are also given. 
The values of $\Gis$ corresponding to the values of $\Gis .f_e$ are these from Fig.~\ref{Fig2}.  }
\label{Fig4}
\end{figure}

A direct check shows that $\Gs_L(\ttt)$ satisfies Langmuir (Frumkin) fractional integral equation (\ref{IBDCAFG}) with error $O({\Gs}^2(\ttt))$. 
In Fig.~\ref{Fig4} both asymptotes, $\Gs_L$ and $\Gs_V$, are compared with numerical results at different values of $\Gis.f_e$. 
The comparisons show that second order asymptotes ($\Gs_V$ and  $\Gs_L$) give good prediction for the adsorption in a broad range ($\Gs\lesssim 3/4$),
while the predictions of the first order asymptotes are satisfactory only at high surfactant concentration ($\Gis.f_e>>1$).
The results in Fig.~\ref{Fig4} confirm also the theoretically derived conclusion (see the discussion below eq. (\ref{RHS}))
that with increasing $\Gis$ ($\Gis.f_e \to 1$, decreasing the initial surfactant concentration)  the adsorption for each of the considered models tends to that of Henry model.  

\begin{figure}[h]%% placement specifier
\centering%% For centre alignment of image.
\includegraphics[width=12cm]{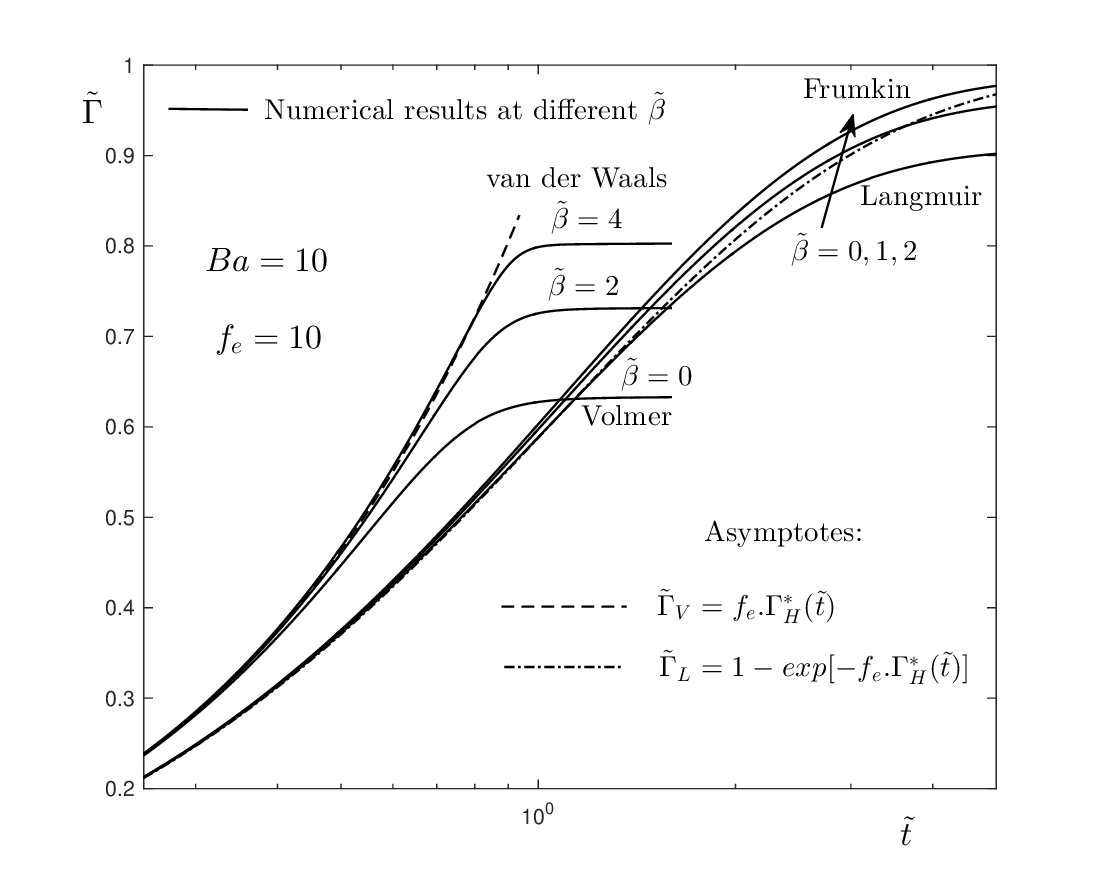}
\caption{Comparison of the second order asymptotes of the adsorption $\tG(\ttt)$ in the case of van der Waals model (dashed line) 
and for Frumkin model (dash-dotted line) 
with numerical results (solid line) at $\Ba=10, \ f_e=10$ and different $\tb$. 
The asymptote $\tG_V(\ttt)$ (dashed line) is the exact solution in case of Henry isotherms. 
} 
\label{Fig5}
\end{figure}

At the end of this subsection results at different values of parameter $\tb$ in the case of van der Waals and Frumkin models are 
presented in Fig.~\ref{Fig5}. 
The results are given in terms of $\tG(\ttt)$ ($\tG(\ttt)=\G(t)/\Gi=\Gs(\ttt)/\Gis$), adsorption scaled with $\Gi$.
This scaling is more natural because $\Gi$ is independent of the models, in contrast to $\G_e$.
The adsorption in terms of $\tG$ is more convenient also for analysis of the models for the surface tension, which are considered in the next subsection. 

The second order asymptotes in terms of $\tG(\ttt)$ are:
    \begin{equation}    
\tG_V(\ttt)=f_e.\Gs_H(\ttt), 
\label{SSTG1}
      \end{equation}
in the case of Volmer (van der Waals) and 
    \begin{equation}    
\tG_L(\ttt)= 1-\exp\left[-f_e.\Gs_H(\ttt)\right], 
\label{SSTG2}
      \end{equation}
in the case of Langmuir and Frumkin models.
Eq. (\ref{SSTG1}) shows that $\tG_H(\ttt)=f_e.\Gs_H(\ttt)=\Gs_H(\ttt).\G_e/\Gi$, which is the exact solution of Henry model, 
is second order asymptote for Volmer and van der Waals models. 
% =f_e.\Gs_H(\ttt)$ ($f_e=\G_e/\Gi$)

The comparisons of the second order asymptotes with numerical results presented so far, including that in Fig~\ref{Fig5}, 
show that the above asymptotes give a good prediction for the adsorption in a wide range.
%(up to about 3/4 of the equivalent value $\G_e$) . 
An important feature of the construction of the second order asymptotes is that their core $\Gs_H(\ttt)$ includes only material parameters. 
The parameters $f_e=K.c_e$ and $\Gis=\Gi/\G_e$ that depend on the initial surfactant concentration $c_e$ and the adsorption isotherm, respectively, 
are multipliers in these asymptotes.

\subsection{Second order asymptotes of the surface pressure}
\label{SSecAsSP}

Second order asymptotes for the surface pressure are derived based on the corresponding asymptotes for the adsorption (\ref{SSTG1}-\ref{SSTG2}).
For this purpose it is convenient to define surface pressure in terms of $\tG(\ttt)=\G(t)/\Gi$:
    \begin{equation}    
\tP(\ttt)=\Pi(t)/(k.T.\Gi)=J(\G(t))/\Gi=J(\tG(\ttt).\Gi)/\Gi  
\label{SPD}
      \end{equation}
Thus plugging the asymptotes $\tG_A(\ttt)$  ($A$ stands for $H$, $V$ or $L$, respectively, see (\ref{SSTG1}-\ref{SSTG2})) 
in the corresponding equation of states given in Table~\ref{Table1}, the asymptotes for $\tP(\ttt)$ are:
\begin{eqnarray}
\tP_H(\ttt) &=& \tG_H(\ttt)= f_e.\Gs_H(\ttt) \label{ASPH}\\
\tP_V(\ttt) &=& \tG_H(\ttt)/[1-\tG_H(\ttt)] -\tb.\tG^2_H(\ttt)/2 \label{ASPV}\\
\tP_L(\ttt) &=& \tG_H(\ttt) -\tb . \left\{1-\exp\left[-\tG_H(\ttt)\right]\right\}^2/2 \label{ASPL}
\end{eqnarray}
for Henry, Volmer (van der Waals) and Langmuir (Frumkin) models, respectively.
Taking into account that $\Gs_H(\ttt)$, given by (\ref{SBDCAH1}), is the exact solution for the adsorption 
in the case of Henry model it follows that $\tP_H(\ttt)$ is the exact solution for the dimensionless surface pressure for this model.

\begin{figure}[h]%% placement specifier
\centering%% For centre alignment of image.
\includegraphics[width=12cm]{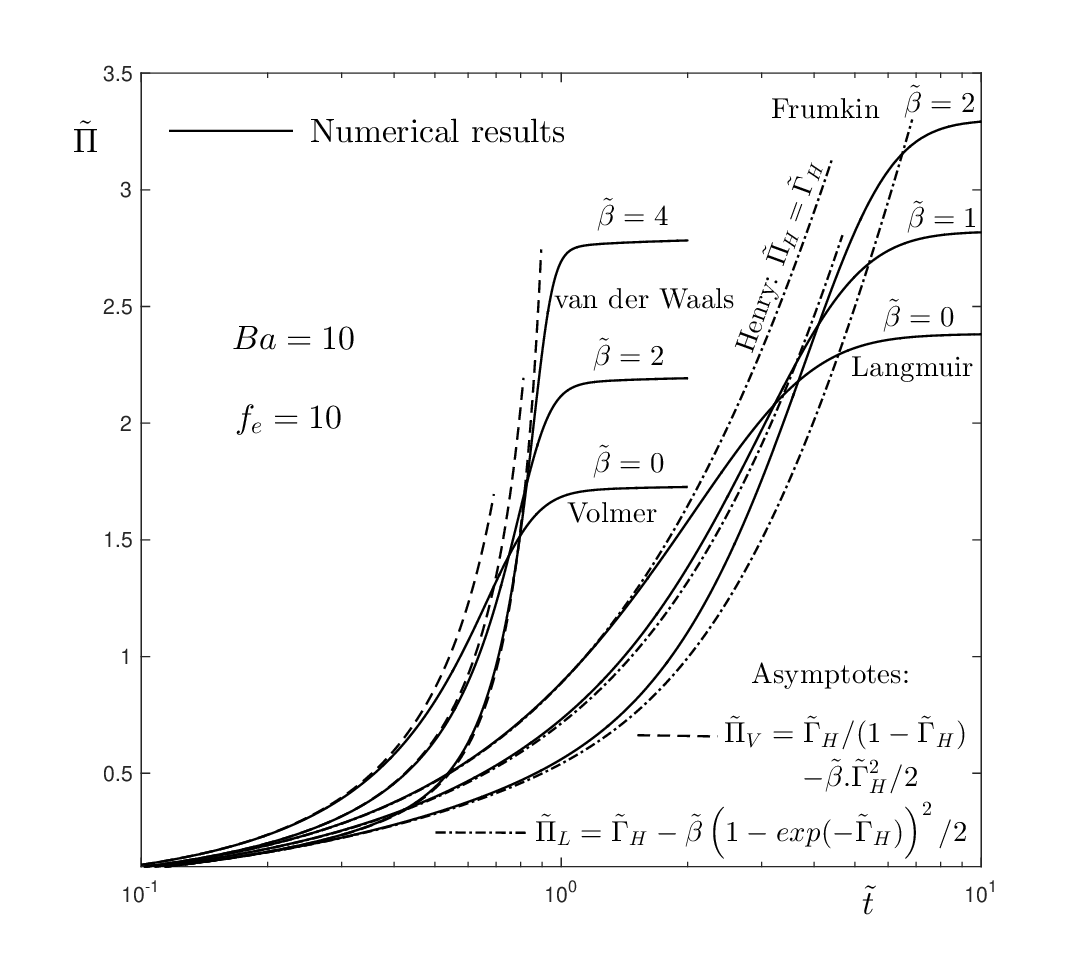}
\caption{Comparison of the asymptotes of the dimensionless surface pressure 
in the case of van der Waals model $\tP_V(\ttt)$ (dashed line) 
and for Frumkin model $\tP_L(\ttt)$ (dash-dotted line) 
with numerical results (solid line) at $\Ba=10, \ f_e=10$ and different values of $\tb$. 
%The asymptote $\tP_L(\ttt)$ (dashed line at $\tb=0$) is the exact solution in case of Henry model $\tP_H$. 
} 
\label{Fig6}
\end{figure}

Because the asymptotes (\ref{SSTG1}-\ref{SSTG2}) for the adsorption are of second order, it follows 
that the asymptotes (\ref{ASPV}-\ref{ASPL}) for the surface pressure are of second order too (with relative error $O({\Gs}^2(\ttt)))$. 
This is confirmed by the comparisons of the asymptotes with numerical results presented in Fig.~\ref{Fig6}.
It follows from (\ref{ASPL}) (see also Fig.~\ref{Fig6}) that the solution of Henry model $\tP_H$ is a second order asymptote for the surface pressure in the case of Langmuir model.

One of the disadvantages of the asymptotes (\ref{ASPV}-\ref{ASPL}) is that they depend on the model under consideration (and the value of $\tb$).
Therefore, it could be useful to have an asymptote that is applicable for all models considered here, even if it is of lower order of accuracy.
Having in mind that $x/(1-x)=x+O(x^2)=1-\exp(-x)$, we see that the following approximation
    \begin{equation}    
\tP_A(\ttt) = \tG_H(\ttt) + O(\tG_H^2(\ttt)) = \tP_H(\ttt) + O(f_e^2.{\Gs_H}^2(\ttt))
\label{AsSPH}
      \end{equation}
is applicable for all considered models ($A$ stands for $H$, $V$ or $L$, respectively). 
Taking into account that $\tP_H(\ttt) = f_e.\Gs_H(\ttt)$ it is clear that the relative error of the above approximation is proportional to $ f_e.\Gs_H(\ttt)$.
Thus it can be concluded that the exact solution of Henry model ($\tP_H(\ttt))=f_e.\Gs_H(\ttt)$) is of second order approximation for Langmuir model (see also Fig.~\ref{Fig6}) and first order for the other three models.

\begin{figure}[h]%% placement specifier
\centering%% For centre alignment of image.
\includegraphics[width=12cm]{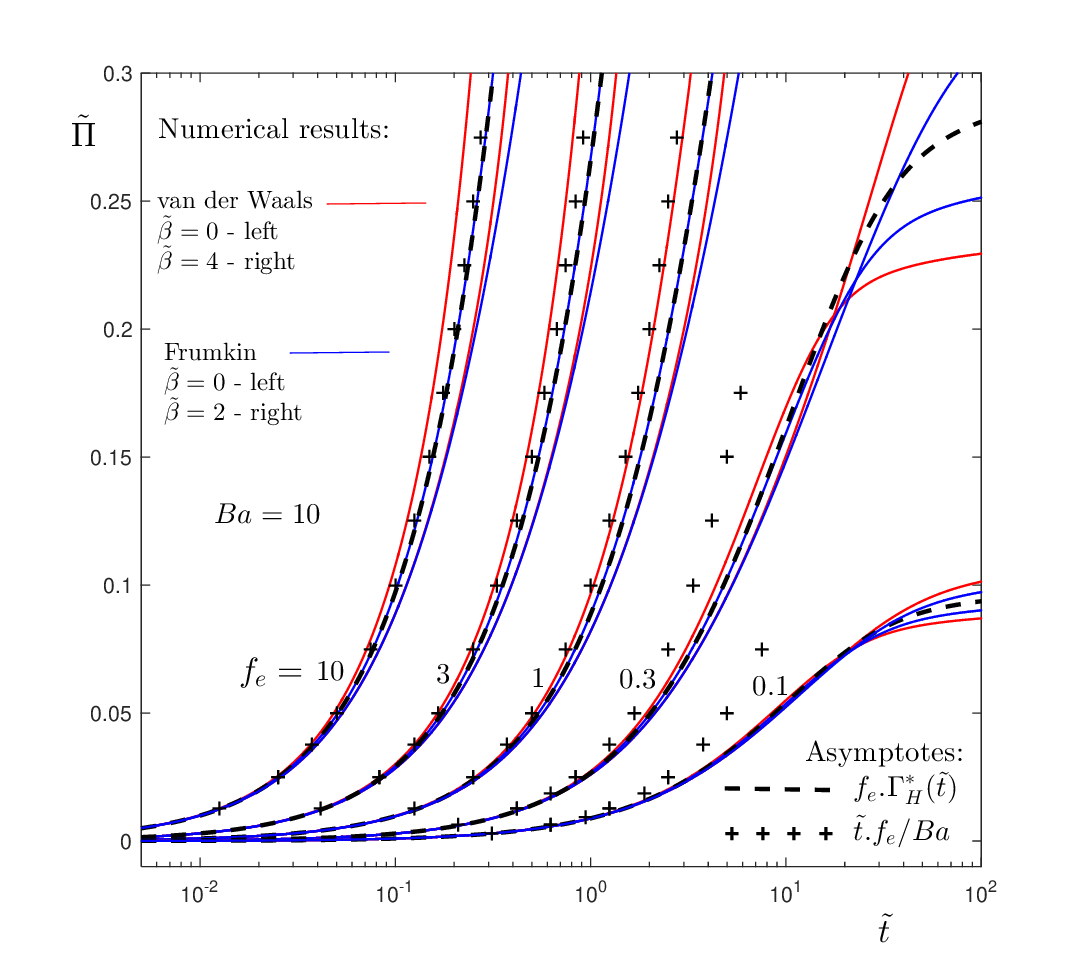}
\caption{Comparison of the first order asymptotes $\tP_H(\ttt)=f_e.\Gs_H(\ttt)$ of the dimensionless surface pressure $\tP(\ttt)$ (dashed lines) 
with numerical results (solid lines) at $\Ba=10$ and different values of $f_e$. 
The marks (+) correspond to the leading term $\ttt/\Ba$ of $\Gs_H(\ttt)$ given by (\ref{SBDCAH3}) (see also (\ref{AsH})).
} 
\label{Fig7}
\end{figure}			

To further investigate this, the results from Fig.~\ref{Fig6}, completed with additional results for different values 
of $f_e$, are shown in Fig.~\ref{Fig7} at the beginning of the process ($\tP<0.3$). 
It is seen that Henry model $\tP_H(\ttt)=f_e.\Gs_H(\ttt)$ (dashed lines)
give good predictions for all models considered here at small values of the surface pressure. 
This is valid independently of the surfactant concentration $c_e$ ($f_e=K.c_e$) for a broad range of the parameter $\tb$ ($\tb\leqslant 4$ for van der Waals and 
$\tb\leqslant 2$ for Frumkin models). 
It is seen also that with decreasing of the initial surfactant concentration ($f_e=K.c_e << 1$) the solutions for the surface pressure of all 
considered models converge to the solution of Henry model (dashed line).
This is in agreement with similar conclusion presented earlier for the adsorption (see Fig.~\ref{Fig4}) and the above estimate, eq. (\ref{AsSPH}). 

The leading term $\ttt.f_e/\Ba=t.c_e.K_a/\Gi$ of the asymptote $\tP(\ttt)=f_e.\Gs_H(\ttt)$, (see (\ref{AsH})) is also shown by marks (+). 
It is worth noting that this also corresponds to the leading term $t^*$ of the solution $1-\exp(-t^*)$ in the limiting case of barrier-controlled adsorption 
% of the solution  
(see eq. (\ref{SBCAH})). 
This explains why the leading term (+) gives good predictions for the evolution of the surface pressure at the beginning of the process where the 
barrier-control is dominant (see also Fig.~\ref{Fig1}). 
It is also seen that the range, with respect to $\tP$, where these predictions (+) are satisfactory 
increases with increasing of the surfactant concentration $c_e$ ($f_e=K.c_e$). 
A discussion in the following section shows that in this case ($\Ba=10$) at the highest considered value of $f_e$ ($f_e=10$) 
the adsorption is dominated by barrier control. 

\section{Numerical approach}
\label{SecNM}

In this section the numerical procedure for computer simulation of the adsorption evolution is discussed. 
It is based on a generalization of the classical Adams-Bashforth-Moulton predictor-corrector method for Volterra type equations. 
The method is widely used for numerical approximation of general fractional equations (see \cite{BB2022} and references therein). 
As a starting point the integral models of the adsorption derived in Section~\ref{SecFIM} are used. 
Important characteristics regarding numerical treatment of the integrals of fractional order is that they are of 
convolution type with singular kernels. 
This requires additional efforts for the numerical approximation and the following computer realization. 
That is why a special attention is paid below to simple steps, that appear to be very important for the accuracy and efficiency of the method.

The core of the numerical approximation of the governing equations (\ref{DCAIS}-\ref{DFG}) is the approximation of the integrals.
Here standard quadrature formulas for convolution type integral with kernel $\xi (t)$ are used:
    \begin{eqnarray}    
\int_0^{t_n} \xi(t_n-\tau) u(\G(\tau)) \, d\tau &\approx& \sum_{j=0}^n u(\G(t_j))\int_0^{t_n}\xi(t_n-\tau) \phi^i _{j}(\tau)d\tau \nonumber\\
&= & \sum_{j=0}^n C^i_{j,n} u(\G(t_j)),\label{QF}
      \end{eqnarray}
where the choice of the basic functions $\phi^i_{j}(\tau)$ ($i=1,2$) determines the order of approximation of the integral.
For the first order approximation the product rectangle rule is applied, where $\phi^1 _{j}(\tau)$ is piecewise step function:
$$ 
\phi_{j}^1(\tau)=\left\{\begin{array}{l}
\ds 1, \ \ \ \tau\in (t_{j},t_{j+1}), \\ [12pt] 
0, \ \ \ \  \mbox{otherwise}.
\end{array}\right.
$$ 
The product trapezoidal quadrature formula gives second order approximation. 
In this case the basic functions are piecewise linear:
$$ 
\phi^2_{j}(\tau)=\left\{\begin{array}{l}
\ds \frac{\tau-t_{j-1}}{t_j-t_{j-1}}, \ \ \ \tau\in (t_{j-1},t_j], \\ [12pt] 
\ds \frac{t_{j+1}-\tau}{t_{j+1}-t_j}, \ \ \ \tau\in (t_{j},t_{j+1}), \\ [12pt]
0, \ \ \ \  \mbox{otherwise}.
\end{array}\right.
$$ 

The singularity at $\tau=t_n$ in the last ($j=n$) integral in (\ref{QF}) is integrable and all integrals have analytical solution 
for the considered in the study kernels $\xi(t)$ and quadratures ($\phi^i_j(\tau)$), see \cite{BB2022}. 
Thus special procedure to deal with the singularity, as in \cite{N2022}, is not necessary here.
Having the coefficients $C^i_{j,n}$, the fractional integral model (\ref{DCAIS}-\ref{IBDCAFG}) is discretized, leading to a relation between $[\G(t_j)]_{j=0}^n$. 
Important characteristics of the relation is that in general it is nonlinear. 
Also, using first order quadrature ($\phi^1_{j}$) the relation is explicit with respect to $\G(t_n)$, while using in the second order ($\phi^2_{j}$) it is implicit.

These features of the quadratures are combined in a predictor-corrector scheme for solving the discretized fractional integral equations. 
Briefly, it combines the explicit predictor step based on rectangle quadrature formula that predicts the solution at the next time step $\G_P(t_n)$. 
This solution is used in the following implicit correction step, based on trapezoidal quadratures, which is of second order accuracy. 
Thus, the advantages of the scheme is second order accuracy and good stability, without necessity of additional iterations typical for the implicit schemes  
(see for instance \cite{LFB2003,CCE2010,LRP1996}).
More information about the numerical scheme can be found in \cite{BB2022}.

The main challenge in the computer implementation of this numerical procedure is due to the convolution type of the fractional integrals.
Numerical approximation of such integrals require $O(n^2)$ operations and information for the solution and the coefficients at all previous time instances 
($\G(t_j), \ C_{j,n}^i, \ 0\leqslant j \leqslant n$). 
At the same time the requirement of accuracy (time steps of order $10^{-5}$) at the beginning and the wide time interval of the process ($t\approx 10^4$) lead to enormous values of $n$ (of order $10^9$). 
To overcome this problem meshes with different time steps in different time intervals are used, 
which appears to be a very efficient approach that significantly improves the performance of the computations. 
In the present study simulations begin with time steps of order $10^{-5}$ to reach values of order $10^{-1}$ at the end where the gradient of $\G$ tends to $0$.

Another challenge of the simulation of the adsorption is that the evolution of the process is very sensitive to the values of the parameters.
It can be seen from the presented so far results that the gradient of the solution is sensitive to $f_e$ and $\Ba$ in a very broad time interval. 
For instance in the limit of diffusion controlled-adsorption, $\Ba<<1$ and $f_e>>1$, (see Fig.~\ref{Fig2})
the highest gradient of the adsorption is at the beginning ($\ttt$ of order $10^{-4}-10^{-5}$). 
In the other limiting regime, barrier-controlled (see Fig.~\ref{Fig1} at $\Ba=10^3$) the gradient of the adsorption is still significant at $\ttt \approx 10^3$. 
Thus, in order to achieve a good approximation the discretization mesh should also depend on the parameters. 
This dependence, however, is not known in advance and it could happen that a chosen in advance mesh is inappropriate for the considered set of parameters. 

To overcome this problem additional rescaling of time is used. 
The idea is that after the rescaling the dynamics of the adsorption (gradients of $\G^*$) happens in a time interval that is less dependent on the parameters. 
To find a proper (in the above described sense) rescaling let us analyze the leading term in the time series of the asymptotes (\ref{SSOABDC1}-\ref{SSOABDC2}). 
Asymptotic expansion (\ref{AsH}) implies that the leading term in the asymptotes for all models is $ \ttt . \Gis.f_e/\Ba$.
Thus, if the time is additionally rescaled as $\ts = \ttt .\Gis.f_e/\Ba$ the leading term of the asymptotes for $\Gs(\ts)$ would be $\ts$, i.e. parameter independent. 
It appears that this is the nondimensionalization of time by the adsorption time $T_a$ (see eq. (\ref{TransTa})).
Thus the diffusion and the adsorption characteristic times are related by $T_d/T_a = \Gis .f_e/ \Ba$. 
The governing  fractional derivative equation (\ref{DBDCA}) in terms of $\Gs (\ts)$, taking into account (\ref{FOTrans}), is:
    \begin{equation}    
\frac{1}{\Phi(\Gs(\ts).\G_e)}\, \D^1_{\ts} \Gs(\ts) + \frac{1}{\sqrt{\Gis.f_e.\Ba}}\D^{1/2}_{\ts} \Gs(\ts) = 1 -\frac{f(\Gs(\ts).\G_e)}{f_e}
\label{DBDCAs}
      \end{equation}

\begin{figure}[h]%% placement specifier
\centering%% For centre alignment of image.
\includegraphics[width=12cm]{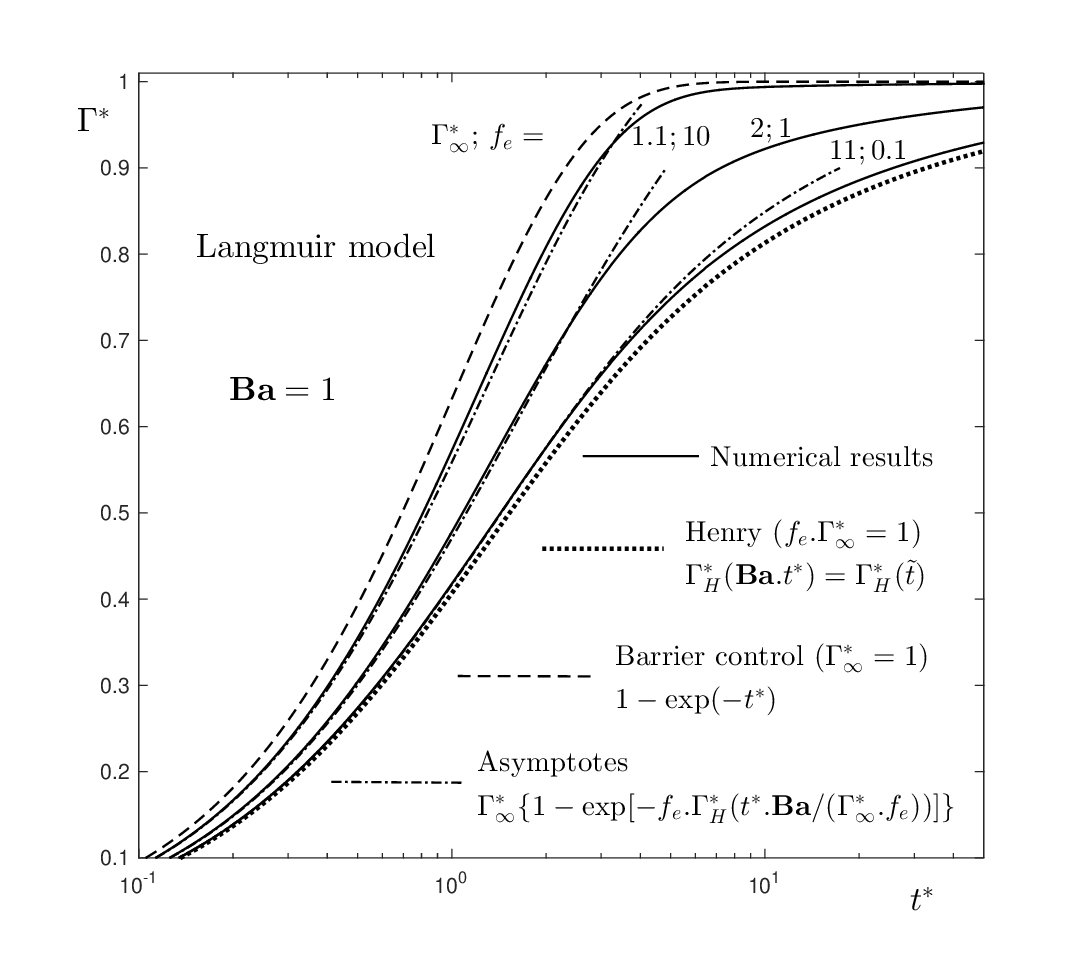}
\caption{Numerical results for Langmuir model of the dimensionless adsorption $\Gs$ as function of the dimensionless time $\ts$ (solid lines). 
The second order asymptotes are given by the dash-dotted lines, as well as the limiting cases of barrier-controlled (dashed line) 
and Henry (dotted line).
} 
\label{Fig8}
\end{figure}			
To demonstrate the effectiveness of the above rescaling, in Fig.~\ref{Fig8} results for $\Gs(\ts)$ are presented in the case of Langmuir kinetic model at $\Ba=1$  and values of $f_e$. 
The results show that independently of the values of $f_e$ the adsorption $\Gs(\ts)$ is between that of Henry model (dotted line) and the
limiting case of barrier-controlled adsorption (dashed line, $1-\exp(-\ts)$). 
This is also valid in the case of Volmer model, with a difference that the upper boundary for the solutions of the adsorption is given by: 
$$ 
\Gs(\ts)=\left\{\begin{array}{l}
\ts, \ \ \ \ts\le 1; \\ [12pt] 
1, \ \ \ \ts \ge 1,
\end{array}\right.
$$ 
which is the limit of barrier-controlled adsorption at $f_e\to \infty$. 

Similar conclusions, that the adsorption $\Gs(\ts)$ is between that predicted by Henry model and the barrier-controlled are valid in the case 
of Frumkin and van der Waals models at $\bt \leq Q$ (for definition of $Q$ see eq. (\ref{TE})). 
This is related with the fact that $\Gis.f_e \geq 1$ for all considered models if $\bt \leq Q$. 

Regarding the numerical simulations important conclusion from the above discussion is that the gradient of the solution at the beginning of the process 
($(d\Gs(\ts)/d\ts)_{\ts=0}=1$) is independent of the parameters of the problem and decreases with time.  
%(at the beginning of the process, $\ts<<1$)the asymptotes $\Gs_L(\ts)$ and $\Gs_V(\ts)$,
%respectively the solution $\Gs(\ts)$, are as $\ts$ at the beginning, $\ts<<1$. 
Also in most of the cases simulations in the time interval $\ts\in [0,1000]$ are sufficient (at $\ts \thickapprox 1000$ the process is already at equilibrium). Simulations for higher $\ts$ are necessary at very small values of $\Gis.f_e.\Ba$ ($\lesssim 10^{-2}$), where the process is already in the diffusion-controlled regime. 
%or relatively large $\bt\gtrsim 3Q$, which is close to the upper limit of $\bt$.
Thus, for the numerical simulations time steps $\Delta\ts\thickapprox 10^{-5}$ at the beginning, to 
$\Delta\ts\thickapprox 0.1$ at the end of the simulations ($\ts\gtrsim 100$), are used in most of the cases. 
The CPU time for a single run in the interval $\ts\in [0,1000]$ is less than 0.1 seconds, using single i7-4710MQ CPU 2.50GHz processor. 
Comparisons with much smaller (10-100 times) time steps indicate that using the above-described meshes the error of the simulations is less than $0.1\%$.

\section{Discussion and conclusions}
\label{SecDC}

In the present study surfactant adsorption under mixed barrier-diffusion control is analyzed by the use of fractional calculus.
This approach simplifies the mathematical model by transforming a system of partial differential equations 
into a single ordinary fractional differential equation (\ref{FBDCA}) for the adsorption.
There are several advantages of this approach, some of which are discussed here.
It leads to a reduction of the dimension of the mathematical model and, respectively, of the number of unknowns.
This facilitates the analysis of the effect of different elements of the model on the adsorption, as well as simplifies the numerical techniques for computer simulations of the evolution of the process.

Having all elements of the model in a single equation allows for a transformation of the variables to a dimensionless form that reduces the number of parameters.
In the kinetic models considered here the dimensionless equation for the adsorption has two dimensionless groups $\Ba =D/(K_a .\Gi .K)$ 
and $f_e = K . c_e$, with a third parameter $\tb = \beta.2.\Gi/(k.T)$ in the cases of Frumkin and van der Waals models.
This further simplifies the numerical simulations and the analysis of the results.

The study is based on the consideration of the integro-differential equation (\ref{BAFD}) as a fractional differential equation (\ref{DBDCA}).
By applying the apparatus of fractional calculus some well-known results are derived directly and new results are obtained. 
An example is a generalization of the Ward-Tordai integral equation to the case of mixed barrier-diffusion controlled adsorption. 
Briefly, the linear combination of derivatives of order $1$ and $1/2$ ($\D^1_{\ttt}$ and $\D^{1/2}_{\ttt}$) in (\ref{DBDCAF1}) 
is considered as a composite fractional derivative $\D^{1,1/2}_{\Ba,\ttt}$.
From fractional calculus theory it is known that the inverse operator to $\D^{1,1/2}_{\Ba,\ttt}$ is the integral operator 
$\J^{1,1/2}_{\Ba,\ttt}$, given in terms of the Mittag-Leffler function $E_{1/2}$.
In this way, the governing equation of the adsorption (\ref{BAFD}) in the case $\Phi (\G)=1$ (see (\ref{DBDCAF1})) 
is transformed to the following integral equation for $\G(t))$:
\begin{eqnarray}
\G(t) &=& \sqrt{D}\int_0^{t}\xi_A (t-\tau)\, \left[c_e-f(\G(\tau))/K\right]\, d\tau, \label{GWT} \\ 
      &=& \sqrt{D}\int_0^{t} A^{-1}\exp \left((t-\tau)/A^2\right)\, \er\left(\sqrt{t-\tau}/A\right)\, \left[c_e-f(\G(\tau))/K\right]\, d\tau, \nonumber
\end{eqnarray}
where the kernel $\xi_A(t)$ is given in (\ref{TTK}) and $A=\sqrt{D}/K_a$. 
This is in fact equation (\ref{IDCAF1}) written in the original dimensional variables $\G$ and $t$.
%As it was discussed in regard to eq. (\ref{IDCAF1}), see (\ref{KernBa}) 

Since in the limiting case $\sqrt{D}/K_a=0$ the kernel is (see (\ref{as2}))
$$\left. \xi^{DC}_0 (t) =\xi_A (t) \right|_{A \to 0} = \frac{1}{\sqrt{\pi t}},$$
equation (\ref{GWT}) reduces to Ward-Tordai equation (\ref{WTFI}), i.e. the case of diffusion-controlled adsorption.
Thus, the integral equation (\ref{GWT}) is a generalization of Ward-Tordai integral equation to the case of mixed barrier-diffusion controlled adsorption.

The other limiting case, of barrier-controlled adsorption (\ref{BCA}), is achieved at kernel $\xi^{BC}_A (t) = 1/A = K_a/\sqrt{D}$.
As it is known, the process of adsorption is barrier-controlled at the beginning, while at the end ($\G(t)$ approaching $\G_e$) 
it is diffusion-controlled.
This is confirmed by the asymptotic behavior of $\xi_A (t)$ given in eqs. (\ref{as1}-\ref{as2}).
They show that at small $t$ ($t \to 0$) $\xi_A(t)$ tends to $\xi^{BC}_A (t)=1/A$, 
while at large times ($t \to \infty$) $\xi_A(t)$ approaches $\xi^{DC}_0 (t)=1/\sqrt{\pi  t}$.
To get a qualitative evaluation of this behavior let us consider the relative difference between the kernel $\xi_A (t)$ 
and those in the limiting cases ($\xi^{BC}_A (t)$ and $\xi^{DC}_0 (t))$.
From (\ref{as1}-\ref{as2}) it follows that:
\begin{eqnarray}
Er^{BC}(t)=\frac{|\xi^{BC}_A (t) - \xi_A(t)|}{\xi^{BC}_A (t)}&=& \frac{2\sqrt{t}}{\sqrt{\pi} A} + O\left(\frac{t}{A}\right), \ \ \ \ \frac{t}{A}<<1, \label{LimBC}\\ 
Er^{DC}(t)=\frac{|\xi^{DC}_0 (t) - \xi_A(t)|}{\xi^{DC}_0 (t)}&=& \frac{A^2}{2 t}+O\left(\frac{A^4}{t^{3/2}}\right), \ \ \ \ \frac{A^4}{t^{3/2}}<<1. \label{LimDC} 
\end{eqnarray}

%Having the above estimates 
Let $T_{BC}$ be the time such that for $t<T_{BC}$ the adsorption $\G(t)$ is approximated by the solution of the model 
of barrier-controlled adsorption (\ref{BCA}). 
Respectively, let $T_{DC}$ be such that at $t>T_{DC}$ the adsorption $\G(t)$ is approximated by that predicted by Ward-Tordai integral equation 
(diffusion-controlled regime).
It follows from (\ref{as1}-\ref{as2}) that $T_{BC}<<D/K_a^2<<T_{DC}$ and the whole adsorption process can be divided in three regimes:
$t \lesssim T_{BC}$ - barrier-controlled; $T_{BC}\lesssim t \lesssim T_{DC}$ - mixed adsorption and $t \gtrsim T_{DC}$ of diffusion-controlled regime.
Now, based on (\ref{LimBC}-\ref{LimDC}), estimates for $T_{BC}$ and $T_{DC}$ can be obtained. 
Thus, at a chosen relative tolerance $Er$ the boundaries of the regimes are given as: 
$$T_{BC}\approx \frac{Er^2\sqrt{\pi}}{2} \frac{D}{K_a^2},\ \ \ \ \ \ \ T_{DC}\approx \frac{1}{2.Er} \frac{D}{K_a^2}.$$
This is also confirmed by the results presented in Fig.~\ref{Fig1} in the case of Henry model at different values of $\Ba$ 
(in this case $\ttt/\Ba= t.K_a^2/D$).

Let us note that the above conclusions are valid in the case $\Phi(\G)=1$. In the other cases (of Langmuir and Frumkin adsorption isotherms) this approach is not applicable and further research is needed. The discussion in the rest of this section concerns all of the considered models.

Another important question is: At what values of the parameters the whole process of adsorption can be considered as barrier-controlled or diffusion-controlled? 
%For the limit of of barrier-controlled adsorption the answer is given based 
First, let us note that at $\Ba. \Gis . f_e \to \infty$ the second term  in eq. (\ref{DBDCAs}) can be neglected and the equation converges to that of barrier-control regime (\ref{BCAD}).
Essential for this analysis is that the right-hand side of the equations is of order $1$ at the beginning.
To give an answer in the other limiting case of diffusion-controlled adsorption, let us rescale the time as $\ttt=\ttt^*.(\Ba. \Gis . f_e)$. 
Equation (\ref{DBDCAs}) in terms of $\ttt^*$ reads:
    \begin{equation}    
\frac{\Ba^*}{\Phi(\tG^*(\ttt^*).\G_e)}\, \D^1_{\ttt^*} \tG^*(\ttt^*) + \D^{1/2}_{\ttt^*} \tG^*(\ttt^*) = 1 -\frac{f(\tG^*(\ttt^*).\G_e)}{f_e},
\label{DDDCAs}
      \end{equation}
where $\Ba^*=\Ba. \Gis . f_e $. 
%where $\tG^*(\ttt^*)=\Gs(\ts)$.
This implies that at $\Ba. \Gis . f_e \to 0$ equation (\ref{DDDCAs}) converges to that in the case of diffusion-controlled adsorption - 
equation (\ref{DCADE}), written in terms of $\tG^*(\ttt^*)$.

Therefore, the boundaries between different regimes of adsorption (diffusion-controlled, mixed and barrier-controlled) 
are determined by the values of the dimensionless group $\Ba^*$.
A quantitative estimation of the boundaries can be obtained based on numerical results. 
They indicate that at $10^{-2} \lesssim \Ba^*=D.c_e/(K_a.\G_e) \lesssim 10^2$ the adsorption is mixed barrier-diffusion controlled.
At smaller values of $\Ba^*$ the process can be considered as diffusion-controlled and at higher as barrier-controlled.
In the case of Henry model, where $\Ba^*=\Ba=D/(K_a.\Gi .K)$, this is seen in Fig.~\ref{Fig1}.
Numerical check shows that these estimates are not influenced by the values of $\Phi(\G(t))$, which can be significant at $\G(t)$ close to $\Gi$ 
in the cases of Langmuir and Frumkin models.

A significant contribution of the presented fractional calculus approach is the derivation of second order asymptotes (at small surface coverage) 
for the adsorption and the surface tension (see Section~\ref{SecAs}).
The asymptotes for the adsorption are based on an approximation of the corresponding function $f(\G(t))$ in the isotherms by that of Henry model 
(see eqs. (\ref{TE}-\ref{RHS}).
In original dimensional variables the asymptotes 
for the adsorption $\G_V(t)$ and $\G_L(t)$ in the cases of Volmer (van der Waals) and Langmuir (Frumkin) models, respectively, are (see (\ref{SSOABDC1}-\ref{SSTG2})):
\begin{eqnarray}
%\G_H(t)&=&\Gi.K.c_e.\Gs_H(\ttt)=\G_e.\Gs_H(\ttt),\nonumber\\ 
\G_V(t)&=&\Gi.K.c_e.\Gs_H(\ttt),\nonumber\\ 
\G_L(t)&=& \Gi\left\{1-\exp\left[-K.c_e.\Gs_H(\ttt)\right]\right\},\nonumber
\end{eqnarray}
Plugging them in the corresponding equations of state, given in Table \ref{Table1}, the asymptotes for the surface tension $\sigma(t)$ 
for Henry, Volmer (van der Waals) and Langmuir (Frumkin) models are:
\begin{eqnarray}
\sigma_H(t)&=& \sigma_0 - k.T.\G_H(t),\nonumber\\ 
\sigma_V(t)&=& \sigma_0 - k.T.\frac{\Gi.\G_H(t)}{\Gi-\G_H(t)} + \beta.\G_H^2(t),\nonumber\\ 
\sigma_L(t)&=&\sigma_0 - k.T.\G_H(t)  + \beta.\Gi^2\left\{1-\exp\left[-\G_H(t)/\Gi\right]\right\}^2,\nonumber
\end{eqnarray}
where $\G_H(t)=\Gi.K.c_e.\Gs_H(\ttt)$ is the exact solution for the adsorption in the case of Henry model. 
The above asymptotes for the surface tension can also be derived by the corresponding asymptotes for the surface pressure (\ref{ASPH}-\ref{ASPL}).

%($\Gi.K.c_e=\G_e$) $\G_H(t)=\G_e.\Gs_H(\ttt)$ is the exact solution for the adsorption, 
%where is the solution in dimensionless variables, .
The core of the asymptotes presented above is the exact solution $\Gs_H(\ttt)$ of the dimensionless Henry model (\ref{SBDCAH2}-\ref{SBDCAH1}).
Important characteristics of this solution is that it depends on only one parameter - the dimensionless group $\Ba$, 
which includes only material parameters ($\Ba=D/(K_a.\Gi.K))$.
The dimensionless time $\ttt$ also depends on material parameters only ($\ttt= t. D/(K.\Gi)^2$).
The initial surfactant concentration $c_e$, which usually is a parameter of the experiment, is a multiplier in the above asymptotes. 
Thus, at a given system of fluids and surfactant ($D,\ K,\ K_a,\ \Gi$), $\Gs_H(\ttt)$ can be calculated in advance for the corresponding value of $\Ba$. 
Then $\Gs_H(\ttt)$ can be used for generation of the asymptotes for any of the considered models at any initial surfactant concentration $c_e$ 
by simple mathematical operations.

The presented asymptotes are validated by comparisons with numerical results at different values of parameters $\Ba,\ f_e$ and $\tb$,
which show that the second order asymptotes give good predictions in a broad interval, up to about $2/3$ of the whole range.
It is interesting to mention that when the adsorption is considered, the exact solution of Henry model $\G_H(t)$ 
is a second order asymptote for Volmer and van der Waals models ($\G_L(t)=\G_H(t)$).
For the surface tension the exact solution of Henry model $\sigma_H(t)$ is second order asymptote for Langmuir model 
and an universal (for all models considered here) first order asymptote.

The approach presented in the study can be straightforwardly applied to other adsorption models, such as adsorption on spherical or expanding interface, at anomalous diffusion, in liquid-liquid systems, at surfactant concentration above CMC.

\section*{Acknowledgments}

The work was partially supported by the Centre of Excellence in Informatics and ICT under the Grant No BG16RFPR002-1.014-0018-C01, financed by the Research,Innovation and Digitalization for Smart Transformation Programme 2021-2027 and co-financed by the European Union.}

%\conflictsofinterest{The authors declare that they have no known competing financial interests or personal relationships that could have appeared to influence the %work reported in this paper.} 

\bibliographystyle{elsarticle-num} 
\bibliography{PaperJCIS2025}



%% For numbered reference style

\end{document}